\documentclass[a4paper,11pt]{amsart}
\usepackage[utf8]{inputenc}

\usepackage{amsmath,amsfonts,amssymb,amsthm,amscd,latexsym,mathrsfs, mathabx, amsbsy,epsfig}

\usepackage{mathtools}
\usepackage[usenames, dvipsnames]{color}

\usepackage{enumitem}
\usepackage{fancyhdr}
\usepackage{stmaryrd, geometry, xypic, tikz-cd, verbatim}

\theoremstyle{definition}
\newtheorem{definition}{Definition}[section]
\newtheorem{proposition}[definition]{Proposition}
\newtheorem{lemma}[definition]{Lemma}
\newtheorem{theorem}[definition]{Theorem}
\newtheorem{corollary}[definition]{Corollary}

\newtheorem{remark}[definition]{Remark}

\newtheorem{example}[definition]{Example}
\newtheorem{claim}[definition]{Claim}

\newcommand{\N}{\mathbb{N}}
\newcommand{\R}{\mathbb{R}}
\newcommand{\Q}{\mathbb{Q}}
\newcommand{\C}{\mathbb{C}}

\newcommand{\M}{\mathfrak{M}}
\newcommand{\Int}{\text{Int}}
\newcommand{\rad}{\text{rad}}
\newcommand{\sbullet}{{\scriptscriptstyle\bullet}}

\DeclareMathOperator{\spec}{Spec}

\DeclareMathOperator{\sym}{Sym}
\DeclareMathOperator{\dist}{dist}

\newcommand{\ndot}{\mathord{\cdot}}
\DeclarePairedDelimiter{\norm}{\lVert}{\rVert}
\DeclarePairedDelimiter{\abs}{\lvert}{\rvert}
\DeclarePairedDelimiter{\algnorm}{\vvvert}{\vvvert}
\DeclarePairedDelimiter{\Card}{\{0,\dots,}{\}}
\DeclarePairedDelimiter{\card}{\{1,\dots,}{\}}
\DeclarePairedDelimiter{\parenth}{\{}{\}}

\begin{document}
\title{Non-Archimedean metric extension for semipositive line bundles}

\author{Yanbo FANG}
\address{IMJ-PRG, Universit\'e Paris Diderot\\75205 PARIS Cedex 13\\France}
\email{yanbo.fang@imj-prg.fr}

\begin{abstract}
    For a projective variety $X$ defined over a non-Archimedean complete non-trivially valued field $k$, and a semipositive metrized line bundle $(L, \phi)$ over it, we establish a metric extension result for sections of $L^{\otimes n}$ from a sub-variety $Y$ to $X$. We form normed section algebras from $(L, \phi)$ and study their Berkovich spectra. To compare the supremum algebra norm $\algnorm{\ndot}_{\phi_Y}$ and the quotient algebra norm $\algnorm{\ndot}_{\phi,X|Y}$ on the restricted section algebra $V_{\sbullet}(L_{X|Y})$, two different methods are used: one exploits the holomorphic convexity of the spectrum, following an argument of Grauert; another relies on finiteness properties of affinoid algebra norms.
\end{abstract}
\maketitle
\tableofcontents

\section{Introduction}
Let $k$ be a field, $X$ be a projective scheme over $\spec k$ and $L$ be an ample invertible $\mathscr{O}_X$-module. Let $Y$ be a closed subscheme of $X$ and $\mathscr{I}_Y$ be the corresponding ideal sheaf. Recall that Serre's vanishing theorem (see for example \cite[Th\'eor\`eme III.2.2.1]{EGA}) asserts that there exists an integer $n_Y$ such that \[H^1(X,\mathscr{I}_Y\otimes L^{\otimes n})=\parenth{0}\] for any integer $n\geqslant n_Y$. Therefore the exact sequence of coherent $\mathscr{O}_X$-modules
\[\xymatrix{\boldsymbol{0}\ar[r]&\mathscr{I}_Y\otimes L^{\otimes n}\ar[r]&L^{\otimes n}\ar[r]&(\mathscr{O}_X/\mathscr{I}_Y)\otimes L^{\otimes n}\ar[r]&\boldsymbol{0}}\]
induces a surjective $k$-linear map from $H^0(X,L^{\otimes n})$ to $H^0(X,(\mathscr{O}_X/\mathscr{I}_Y)\otimes L^{\otimes n})$, which coincides with the restriction map $H^0(X,L^{\otimes n})\rightarrow H^0(Y,L|_Y^{\otimes n})$ if we identify $H^0(Y,L|_Y^{\otimes n})$ with $H^0(X,(\mathscr{O}_X/\mathscr{I}_Y)\otimes L^{\otimes n})$. In other words, for any integer $n\geqslant n_Y$, any global section of $L|_Y^{\otimes n}$ extends to a global section of $L^{\otimes n}$. We denote by $H^{0}(Y, L_{X|Y}^{\otimes n})$ the image vector space of restriction map
\[H^{0}(Y, L_{X|Y}^{\otimes n}):=\mathrm{Image}(H^0(X,L^{\otimes n})\longrightarrow H^0(Y,L|_Y^{\otimes n})).\]

Assume that $k$ is equipped with a complete absolute value $\abs{\ndot}$ and that $L$ is equipped with a continuous metric $\phi=(\abs{\ndot}_{\phi}(x))_{x\in X^{\mathrm{an}}}$, which induces by tensor power a continuous metric $n\phi$ on each $L^{\otimes n}$, $n\in\mathbb N$. Suppose in addition that the schemes $X$ and $Y$ are integral. Then the space of global sections $H^0(X,L^{\otimes n})$ is naturally equipped with a supremum norm $\norm{\ndot}_{n \phi}$ associated with $n \phi$, defined as follows 
\[\forall\,s\in H^0(X,L^{\otimes n}),\quad\norm{s}_{n\phi}:=\sup_{x\in X^{\mathrm{an}}}|s(x)|_{n\phi}(x).\]
We denote by $\phi|_Y$ the restriction of the metric $\phi$ on $L|_Y$. A supremum norm $\norm{\ndot}_{n\phi|_Y}$ on $H^0(Y,L|_Y^{\otimes n})$ is defined in a similar way. The metric extension problem compares the norm $\norm{\ndot}_{n\phi|_Y}$ to the quotient norm (denoted by $\norm{\ndot}_{n\phi, X|Y}$) of $\norm{\ndot}_{n\phi}$ induced by the restriction map $H^0(X,L^{\otimes n})\rightarrow H^0(X,L|_Y^{\otimes n})$ with $n\in\mathbb N$, $n\geqslant n_Y$. Note that by definition we always have $\norm{\ndot}_{n\phi, X|Y}\geqslant \norm{\ndot}_{n\phi|_Y}$ on $H^{0}(Y, L_{X|Y}^{\otimes n})$. Therefore the metric extension problem can be interpreted as finding a uniform upper bound for 
\begin{equation}\label{Equ: uniform upper bound}\inf_{\begin{subarray}{c}s\in H^0(X,L^{\otimes n})\\
s|_Y=t
\end{subarray}}\frac{\norm{s}_{n\phi}}{\norm{t}_{n\phi|_Y}},\quad t\in H^0(Y,L|_Y^{\otimes n})\setminus\{0\}.\end{equation}

In the complex analytic setting, namely when $(k,\abs{\ndot})$ is $\mathbb C$ equipped with the usual absolute value, the metric extension problem has been studied by different authors using various approaches. Assume that the metric $\phi$ is strictly positive (namely, for any local section $s$ of $L$ over an open subscheme $U$ of $X$,  the function $(x\in U^{\mathrm{an}})\mapsto \log|s(x)|_{\phi}$ is strongly plurisubharmonic). In the case where $X$ is smooth, by Gromov's theorem we can compare the sup norm $\norm{\ndot}_{n\phi}$ to the $L^2$ norm $\norm{\ndot}_{n\phi,L^2}$ defined as 
\[\forall\,s\in H^0(X,L^{\otimes n}),\quad\norm{s}_{n\phi,L^2}:=\bigg(\int_{X^{\mathrm{an}}}|s(x)|_{n\phi}(x)\,\mathrm{d}V\bigg)^{\frac 12},\]
where $\mathrm{d}V$ is a probability measure on $X^{\mathrm{an}}$ which is locally equivalent with Lebesgue measure with a smooth Radon-Nikodym density. Therefore, in the case where $X$ and $Y$ are both smooth, we can apply the Andreotti-Vesentini-H\"omander's $L^2$ technique or $L^2$-extension theorems of Ohsawa-Takegoshi type \cite{OT} and get an inequality (see for example \cite{Tia} and \cite{Man}) 
\begin{equation}\label{Equ: Tian bound}\norm{\ndot}_{n\phi|_Y}\geqslant C'(\phi,Y, X)n^{-d}\norm{\ndot}_{n\phi, X|Y},\qquad n\geqslant n_Y,\end{equation}
where $C'(\phi,Y, X)$ is a positive constant. Alternatively, one can apply Grauert's argument of pseudo-convexity of the (open) dual unit disc bundle of $(L, \phi)$ to produce, for any $\epsilon>0$, a slightly weaker inequality of the form
\begin{equation}\label{Equ: Bost-Randriam bound}\norm{\ndot}_{n\phi|_Y}\geqslant C_\epsilon(\phi,Y, X)\mathrm{e}^{-\epsilon n}\norm{\ndot}_{n\phi, X|Y},\qquad n\geqslant n_Y,\end{equation}
where $C_\epsilon(\phi,Y, X)$ are is a positive constant depending on $\epsilon$. We refer to \cite{Bos} and \cite{Ran} for more details.

Note that, to obtain an estimate of the form \eqref{Equ: Bost-Randriam bound} for any closed point $Y$, the semipositivity of the metric $\phi$ is actually necessary. In fact, for sufficiently positive integer $n$, the ample invertible $\mathscr{O}_X$-module determines a closed embedding $\iota_n:X\rightarrow\mathbb P(H^0(X,L^{\otimes n}))$. The norm $\norm{\ndot}_{n\phi}$ on $H^0(X,L^{\otimes n})$ defines a Fubini-Study metric on the universal invertible sheaf of the projective space $H^0(X,L^{\otimes n})$. Denote by $\phi_n$ the continuous metric on $L$ the $n$-th tensor power of which identifies with the restriction of the Fubini-Study metric. Then the estimate \eqref{Equ: Bost-Randriam bound} in the case where $Y$ is a single closed point $\{x\}$ implies that the sequence of norms $\abs{\ndot}_{\phi_n}(x)$ converges to $\abs{\ndot}_{\phi}(x)$. Moreover, it can be shown the subsequence $(|\ndot|_{\phi_{2^m}}(x))_{m\in\mathbb N}$ is decreasing. Therefore, the semipositivity of Fubini-Study metrics implies that of the metric $\phi$. 

The problem of metric extension has various applications, not only in complex analytic geometry, but also in Arakelov geometry. It is a key ingredient in the proof of the arithmetic Hilbert-Samuel theorem, see \cite{AB}. It has also been applied in the proof of Nakai-Moishezon criterion of arithmetic ampleness, cf. \cite{Zhang}, see also \cite{Mor}. 

From the adelic point of view of Arakelov geometry, one can replace the integral models of arithmetic objects by a family of Berkovich analytic objects (possibly equipped with metrics) parametrised by the set of finite places of a number field. The advantage of the adelic approach consists in treating all places of a number field in a uniform way. This motivates the research of the non-Archimedean analogue of notions and results in complex analytic geometry. It turns out that many usual analytic tools (such as $L^2$ estimates) do not work well in the non-Archimedean setting, and often new ideas are needed to develop the non-Archimedean analogue of complex analytic geometry and to unify the arguments in both settings.

In this article we undertake a study of the metric extension problem in (\ref{Equ: uniform upper bound}) in the non-Archimedean analytic setting. Various notions of semipositivity of a metric have been proposed (see \S3.1 \ref{Item: FS metric and FS envelop metric}), we adopt the one as being a uniform limit of Fubini-Study metrics. We establish the following result (see Theorem \ref{quantitative extension geometric}, Theorem \ref{quantitative extension algebraic}), which improves considerably the metric extension theorem of \cite{CM}.
\begin{theorem}\label{quantitative extension}
    Let $\phi$ be an asymptotic Fubini-Study metric on $L$. Then for any $\epsilon>0$, there exists $n_Y\in\N$ such that, for any $n\geq n_Y$ and any $t_n\in H^{0}(Y, L|_Y^{\otimes n})$, there exits $s_n\in V_n(L)$ such that $s_n|_Y=t_n$ and 
    \[\norm{s_n}_{n\phi}\leq \mathrm{e}^{n\epsilon}\cdot\norm{t_n}_{n\phi|_Y}\]
\end{theorem}
We mimic Grauert's argument as in \cite{Bos} and in \cite{Ran}, though some parts of the argument require adaptations. On the restricted section algebra, namely
\[V_{\sbullet}(L_{X|Y}):=\bigoplus_{n\in \N}H^{0}(Y, L_{X|Y}^{\otimes n})=\mathrm{Image}(V_{\sbullet}(L)\longrightarrow V_{\sbullet}(L|_Y)),\]
we gather all norms on graded pieces together, to have two algebra norms defined as follows
\[\forall \underline{t}=(t_n)_{n\in \N}\in V_{\sbullet}(L_{X|Y}), \ \algnorm{\underline{t}}_{\phi|_Y}:=\sup_{n\in \N}\norm{t_n}_{n\phi|_Y}, \ \algnorm{\underline{t}}_{\phi, X|Y}:=\sup_{n\in \N}\norm{t_n}_{n\phi, X|Y},\]
respectively. After completion with respect to these algebra norms, we get two commutative Banach algebras $\widehat{V}_{\sbullet}(L_{X|Y}, \phi|_{Y})$ and $\widehat{V}_{\sbullet}(L_{X|Y}, \phi_{X|Y})$. The geometry of Berkovich spectra of these Banach algebras are studied in \S3. They are intimately related to the dual unit disc bundle of $L$ with respect to the Fubini-Study envelop metric $\mathcal{P}(\phi)$. It turns out that when $\phi$ is semipositive (more precisely, when $\phi$ is asymptotic Fubini-Study), the spectral seminorm of $\algnorm{\ndot}_{\phi, X|Y}$ is $\algnorm{\ndot}_{\phi|_Y}$, and the two spectra coincide. This already permits us to reprove a statement of metric extension of Chen-Moriwaki in the non-trivial valuation case (\cite[Theorem 0.1]{CM}, see Theorem \ref{weak QE}). Note that since we equip $V_{\sbullet}(L_{X|Y})$ with Banach algebra norms instead of Fréchet algebra seminorms as used in \cite{Bos}\cite{Ran}, the spectrum gives rise to closed dual unit disc bundle $\overline{\mathbb{D}}^{\vee}(L|_Y, \mathcal{P}(\phi)|_Y)$ rather than the open dual disc bundle.

There are two independent methods to compare these two algebra norms. The geometric method in \S4 exploits the holomorphic convexity of the aforementioned spectrum, and uses holomorphic functional calculus in non-Archimedean commutative Banach algebras, to construct a (continuous) Banach algebra homomorphism from $\widehat{V}_{\sbullet}(L_{X|Y}, \phi_{X|Y})$ to some perturbed version of $\widehat{V}_{\sbullet}(L_{X|Y}, \phi|_{Y})$. This construction is analogous to the use of Grauert's vanishing theorem along with pseuodo-convexity to get continuous map between Fréchet seminormed coherent analytic sheaves in the $\C$-analytic setting. In \S5 we present an algebraic method with an extra assumption that $(k, \abs{\ndot})$ is discretely valued. The analytic convexity enters only in the semi-positivity of the continuous metric, equivalently the uniform approximation by Fubini-Study metrics. For $\phi$ being a Fubini-Study metric, we show directly by a delicate calculation that the two Banach algebra norms are affinoid algebra norms, which possesses strong finiteness property to give a uniform upper bound. The exact calculation depends heavily on the existence of a non-Archimedean orthogonal basis. We are unaware of any analogue of this affinoid algebra technique in the $\C$-analytic setting, but it seems plausible to compare this with the use of Ohsawa-Takegoshi $L^2$ extension theorem (see Remark \ref{affinoid estimate vs OT L2}).

Note that our proof gives a sub-exponential upper bound as in (\ref{Equ: Bost-Randriam bound}). Whether the polynomial bound of (\ref{Equ: Tian bound}) is obtainable in the non-Archimedean setting remains unexplored. It seems to us that the commutative Banach algebra techniques that we used here are insufficient to ameliorate the bound. We would like also to metion the recent work \cite{MP} which gives the hope of carrying out directly Grauert's argument as in \cite{Bos} for the metric extension problem in the non-Archimedean setting.

\section{Reminders on ultrametric functional analysis}

In this section, one recalls some facts about functional analysis concerning normed vector spaces and normed algebras over a complete non-Archimedean valued field, following \cite{Ber}, \cite{BGR}, \cite{FvdP} and \cite{Tem}. Besides the well-known ones, results in \S2.1.2, \S2.2.4, \S2.3.3 and \S2.4 are most relevant to our construction.

Throughout the section, one fixes a field $k$ equipped with a \emph{non-Archimedean} and non-trivial absolute value $\abs{\ndot}$ and we assume that $k$ equipped with the topology defined by the absolute value is complete. Denote by $k^{\circ}$ the valuation ring of $(k, \abs{\ndot})$, by $k^{\circ\circ}$ the maximal ideal of $k^{\circ}$, and by $\widetilde{k}$ the residual field $k^{\circ}/k^{\circ\circ}$. Denote by $H(k, \abs{\ndot})$ the $\Q$-vector subspace of $\R$ generated by the set of numbers $\log\abs{k^{\times}}$, and by $\alpha$ the quotient map of $\Q$-vector spaces $\R\to \R/H(k, \abs{\ndot})$. One says that $n$ real numbers $\parenth{p_1, \dots, p_n}$ are \emph{$\Q$-independent} in $\R/H(k, \abs{\ndot})$ if the vectors $\parenth{\alpha(p_1), \dots, \alpha(p_n)}$ are $\Q$-linearly independent in $\R/H(k, \abs{\ndot})$. Unless specified, all $k$-algebras are supposed to be commutative unitary (with $0\neq 1$) and by convention all homomorphism of $k$-algebras are supposed to preserve the units. 

\subsection{Seminormed vector spaces}
\subsubsection{Basic constructions}
Let $V$ be a vector space over $k$. By \emph{seminorm} on $V$, we refer to a map $\norm{\ndot}:V\rightarrow\mathbb R_{\geq 0}$ such that $\norm{ax}=\abs{a}\cdot\norm{x}$ for any $(a,x)\in k\times V$ and that $\norm{x+y}\leqslant\norm{x}+\norm{y}$ for any $(x,y)\in V\times V$. The couple $(V,\norm{\ndot})$ is called a \emph{seminormed vector space} over $k$. If in addition $\norm{\ndot}$ takes positive values on $V\setminus0$, we say that $\norm{\ndot}$ is a \emph{norm} on $V$ and that $(V,\norm{\ndot})$ is a normed vector space. Denote by $\mathfrak n(\norm{\ndot})$ the inverse image of $\{0\}$ by $\norm{\ndot}$. It can be shown that $\mathfrak n(\norm{\ndot})$ is a closed vector subspace of $V$, called the \emph{null space} of $\norm{\ndot}$. Note that there exists a unique norm on $V/\mathfrak n(\norm{\ndot})$, the composition of which with the projection map $V\rightarrow V/\mathfrak n(\norm{\ndot})$ identifies with the seminorm $\norm{\ndot}$. Call this norm the \emph{induced norm} of $\norm{\ndot}$.

Let $(V,\norm{\ndot})$ be a seminormed vector space over $k$. If the strong triangle inequality holds for the seminorm $\norm{\ndot}$, namely $\norm{x+y}\leqslant\max(\norm{x},\norm{y})$ for any $(x,y)\in V\times V$, we say that the seminorm is \emph{ultrametric}. Note that if a seminorm is ultra-metric, then the inequality becomes an equality whenever $\norm{x}\neq \norm{y}$.

We say that a seminormed (resp. normed) vector space $(V,\norm{\ndot})$ is \emph{complete}, or $\norm{\ndot}$ is a \emph{complete seminorm} (resp. \emph{complete norm}) on $V$, if any Cauchy sequence in $V$ with respect to the seminorm $\norm{\ndot}$ admits a limit. A complete normed vector space over $k$ is called a \emph{Banach space} over $k$. Any finite-dimensional normed space $(V,\norm{\ndot})$ is complete. (\cite[1.2.3 Theorem 2]{Bou})

Let $(V,\norm{\ndot}_V)$ be a seminormed vector space over $k$. Let $\widetilde{V}_c$ be the vector space of all Cauchy sequences in $V$ with respect to $\norm{\ndot}_V$. We define a seminorm $\norm{\ndot}_c$ on $\widetilde{V}_c$ which sends any Cauchy sequence $\{v_i\}_{i\in \N}$ to $\lim_{i\rightarrow +\infty}\norm{v_i}_{V}$. Denote by $V_c$ the quotient vector space $V_c/\mathfrak n(\norm{\ndot}_c)$. Then the vector space $V_c$ equipped with the norm induced by $\norm{\ndot}_c$ forms a Banach space over $k$, called the \emph{separated completion} of $(V, \norm{\ndot})$. Tautologically it can be shown that this Banach space is canonically isomorphic to the completion of $V/\mathfrak n(\norm{\ndot})$ equipped with the quotient norm induced by the seminorm $\norm{\ndot}$.

\begin{definition}
Let $\norm{\ndot}_1$ and $\norm{\ndot}_2$ be seminorms on $V$. We say that $\norm{\ndot}_1$ and $\norm{\ndot}_2$ are \emph{equivalent} if there exist two constants $C_1>0$ and $C_2>0$ such that $C_1\norm{\ndot}_1\leq \norm{\ndot}_2\leq C_2\norm{\ndot}_1$. Note that this condition  holds if and only if the seminorms $\norm{\ndot}_1$ and $\norm{\ndot}_2$ induce the same topology on the vector space $V$ (\cite[Corollaire I.3.3.1]{Bou})(note that the absolute value $\abs{\ndot}$ is supposed to be non-trivial). 
\end{definition}

\begin{definition}
Let $(V,\norm{\ndot})$ be a seminormed vector space over $k$. If $W$ is a vector subspace of $V$, then map $(x\in W)\mapsto \norm{x}$ defines a seminorm on $W$, called the \emph{restriction} of $\norm{\ndot}$ on $W$. If $Q$ is a quotient vector space of $V$ and $\pi:V\rightarrow Q$ is the quotient map, then the map $(q\in Q)\mapsto\inf_{x\in\pi^{-1}(\{q\})}\norm{x}$ defines a seminorm on $Q$, called the \emph{quotient} of $\norm{\ndot}$ on $Q$. 
\end{definition}

\begin{definition}
Let $(V,\norm{\ndot}_V)$ and $(W,\norm{\ndot}_W)$ be seminormed vector spaces over $k$, and $f:V\rightarrow W$ be a $k$-linear map. We say that $f$ is \emph{bounded} if there exists a constant $C>0$ such that $\norm{f(x)}_W\leq C \norm{x}_V$ for any $x\in V$. Note that this condition holds if and only if $f$ is continuous with respect to the topologies on $V$ and $W$ induced by the seminorms $\norm{\ndot}_V$ and $\norm{\ndot}_W$ respectively. We say that $f$ is \emph{admissible} if it is bounded and if on the image of $f$, the quotient seminorm of $\norm{\ndot}_V$ and the restriction of $\norm{\ndot}_W$ are equivalent. 
\end{definition}

We recall below several fundamental results in functional analysis and refer to \cite[Theorem 1.3.3.1, Corollary 1.3.3.1, 1.3.3.2, 1.3.3.5]{Bou}  for more details.

\begin{theorem}\label{fundamental Banach}
Let $(V,\norm{\ndot}_V)$ and $(W,\norm{\ndot}_W)$ be Banach spaces over $k$, and $f:V\rightarrow W$ be a $k$-linear map.
    \begin{enumerate}[label=\rm{(\arabic*)}]
        \item\label{Item: closed graph} The $k$-linear map $f$ is bounded if and only if its graph in $V\times W$ is closed under the product topology.
        \item\label{Item: open mapping} Assume that $f$ is bounded and surjective, then $f$ is an open map.  In particular, the quotient norm of $\norm{\ndot}_V$ on $W$ is equivalent to $\norm{\ndot}_W$.
        \item\label{Item: closed image} Assume that $f$ is bounded and injective, then $f(V)$ is closed in $W$.
    \end{enumerate}
\end{theorem}

\begin{theorem}\label{equivalent norm}
Let $V$ be a vector space over $k$ and $\norm{\ndot}_1$ and $\norm{\ndot}_2$ be complete norms on $V$. If there exists $C>0$ such that $\norm{\ndot}_2\leq C\norm{\ndot}_1$, then the norms $\norm{\ndot}_1$ and $\norm{\ndot}_2$ are equivalent.
\end{theorem}

Using this norm equivalence theorem for Banach spaces over $k$, we have immediately the following
\begin{corollary}
Let $(V,\norm{\ndot}_V)$ and $(W,\norm{\ndot}_W)$ be Banach spaces over $k$, and $f:V\rightarrow W$ be a bounded $k$-linear map with closed image. Then $f$ is admissible.
\end{corollary}

\begin{definition}\label{dual norm}
    Let $(V, \norm{\ndot})$ be a finite-dimensional normed vector space. The \emph{dual norm of} $\norm{\ndot}^{\vee}$ on the dual vector space $V^{\vee}$ is defined by
    \[\forall \ell \in V^{\vee}, \ \norm{\ell}^{\vee}:=\sup_{v\in V\setminus\parenth{0}}\frac{\abs{\ell(v)}}{\norm{v}}.\]
\end{definition}
\begin{remark}
    The norm $\norm{\ndot}^{\vee}$ is ultrametric, and $\norm{\ndot}^{\vee\vee}=\norm{\ndot}$ if and only if $\norm{\ndot}$ is ultrametric. (\cite[Section 2.2.3]{CM})
\end{remark}
\begin{definition}\label{scalar extension of norm}
    Let $(V,\norm{\ndot})$ be a normed vector space. Let $(k', \abs{\ndot}')$ be a complete valued field extension of $(k, \abs{\ndot})$. Set $V_{k'}$ to be $V\otimes_k k'$, which can be identified with $\mathrm{Hom}_k(\mathrm{Hom}_k(V, k), k')$. The norm
    \[\forall v'\in V_{k'},\ \norm{v'}_{k'}:=\sup\Big\{\frac{\abs{(\ell\otimes 1) (v')}'}{\norm{\ell}^{\vee}},\ \ell\in V^{\vee}\Big\}\]
    defined via this identification is called the \emph{scalar extension of} $\norm{\ndot}$. 
\end{definition}
\begin{remark}
    If $\norm{\ndot}$ is ultrametric, then $\norm{\ndot}_{k'}$ is the largest ultrametric norm on $V_{k'}$ extending $\norm{\ndot}$. (\cite[Definition 2.4]{CM})
\end{remark}
\begin{lemma}\label{rk 1 quotient norm}
    Let $f: V\rightarrow W$ be a surjective $k$-linear map of finite-dimensional vector spaces, with $\mathrm{dim}_k W=1$. Let $\norm{\ndot}_V$ be a norm on $V$ and let $\norm{\ndot}_W$ be its quotient norm for $f$. Then the norm $\norm{\ndot}_{W, k'}$ identifies with the quotient norm of $\norm{\ndot}_{V, k'}$ induced by the surjective $k'$-linear map $f\otimes \mathrm{id}_{k'}: V_{k'}\rightarrow W_{k'}$. (\cite[Lemma 2.5]{CM})
\end{lemma}

\subsubsection{Orthogonal basis}

\begin{definition}
Let $(V,\norm{\ndot})$ be a finite-dimensional normed vector space over $k$. A basis $\{e_i\}_{i\in \card{n}}$ of $V$ is called \emph{orthogonal (with respect to $\norm{\ndot}$)} if
\[\forall (c_1,\dots, c_{r})\in k^{n},\quad\Big\|\sum_{i\in \card{n}}c_ie_i\Big\|=\max_{i\in \card{n}}\norm{c_ie_i}\] 
Moreover, it is said to be \emph{orthonormal} if in addition $\norm{e_i}=1$ for all $i\in \card{n}$.
\end{definition}

\begin{lemma}\label{distinct orthogonal}
    Let $(V,\norm{\ndot})$ be a finite-dimensional ultrametrically normed vector space over $k$. If $\parenth{v_i}_{i\in \card{n}}$ is a finite set of elements of $V$ such that $\parenth{\norm{v_i}}_{i\in \card{n}}$ are disctinct in $\R_+$. Then $\norm{ \sum_{i\in \card{n}} v_i}=\max_{i\in \card{n}}\norm{v_i}$.
\end{lemma}
\begin{proof}
    If $n=2$, this is clear from the ultra-metric inequality. For general $n$ an induction argument shows the equality.
\end{proof}

\begin{corollary}\label{independence orthogonal}
    Let $(V,\norm{\ndot}_V)$ be a finite-dimensional ultrametrically normed vector space over $k$. Suppose that $(k, |\cdot|)$ is discretely valued. If $\parenth{e_i}_{i\in \card{n}}$ is a basis of $V$ such that $\parenth{\log \norm{e_i}}_{i\in \card{n}}$ are $\Q$-independent in $\R/H(k, \abs{\ndot})$, then $\parenth{e_i}_{i\in \card{n}}$ is an orthogonal basis.
\end{corollary}
\begin{proof}
    For any $f=(f_1, \dots, f_n)\in (k^{\times})^{n}$, the numbers $\parenth{\abs{f_i}\norm{e_i}}_{i\in \card{r}}$ are distinct, otherwise there exist $i,j\in \card{n}, i\neq j$ such that \[\log\norm{e_i}-\log\norm{e_j}=\log\Big| \frac{f_i}{f_j}\Big| \in \log\abs{k^{\times}}\]
    which contradicts the assumption of $\Q$-independence. Hence 
    \[\Big\|\sum_{i\in \card{n}}f_ie_i\Big\|=\max_{0\leq i\leq n}\norm{f_ie_i}\]
    by Lemma \ref{distinct orthogonal}.
\end{proof}

\begin{proposition}\label{discrete orthogonal vector space}
    Let $(V,\norm{\ndot}_V)$ be a finite-dimensional ultrametrically normed vector space over $k$. Suppose that $(k, \abs{\ndot})$ is discretely valued. Then there exists an orthogonal basis for $(V, \norm{\ndot})$. (\cite[Proposition 2.5]{BMPS})
\end{proposition}

\subsection{Banach algebra}

\subsubsection{Basic constructions}
\begin{definition}
Let $A$ be a $k$-algebra (the unit of which is denoted by $\mathbf{1}$) and $\norm{\ndot}$ be a seminorm on $A$ (viewed as a vector space over $k$). 
\begin{enumerate}[label=\rm{(\arabic*)}]
\item The seminorm $\norm{\ndot}$ is said to be \emph{sub-multiplicative} if for any $(a,b)\in A\times A$ one has $\norm{ab}\leq\norm{a}\cdot\norm{b}$. 
\item The seminorm $\norm{\ndot}$ is called \emph{power-multiplicative} if $\norm{a^n}=\norm{a}^n$ for any $a\in A$ and any $n\in\mathbb N\setminus\{0\}$. 
\item The seminorm $\norm{\ndot}$ is called \emph{multiplicative} if $\norm{ab}=\norm{a}\cdot\norm{b}$ for any $(a,b)\in A^2$.
\end{enumerate}
A \emph{$k$-algebra seminorm} (resp. \emph{$k$-algebra norm}) on $A$ is defined to be a sub-multiplicative seminorm (resp. sub-multiplicative norm) $\norm{\ndot}$ on $A$ such that $\norm{\mathbf{1}}=1$. We denote by $\algnorm{\ndot}$ an algebra seminorm. Any $k$-algebra equipped with a complete $k$-algebra norm is called a \emph{Banach $k$-algebra}. 

We use calligraphic letters to denote Banach algebras and Banach modules (defined below) and use the corresponding capital letters to denote the underlying $k$-algebra or the underlying module of a $k$-algebra. For example, a Banach $k$-algebra $(A,\algnorm{\ndot})$ is denoted by $\mathcal{A}$. If $A'$ is a sub-$k$-algebra of $A$, then the restriction of $\algnorm{\ndot}$ on $A'$ is a $k$-algebra norm. If  this norm is complete, we say that $\mathcal{A}'$ ($A'$ equipped with the restricted norm) is a Banach $k$-sub-algebra of $\mathcal{A}$. Similarly, if $Q$ is a quotient $k$-algebra of $A$, then the quotient of the norm $\algnorm{\ndot}$ on $Q$ is a sub-multiplicative seminorm. If it is a complete norm, we say that $\mathcal{Q}$ ($Q$ equipped with the quotient norm) is a Banach quotient $k$-algebra of $\mathcal{A}$.
\end{definition}

\begin{example}\label{Tate algebra}
Let $\mathcal{A}$ be a Banach $k$-algebra. The Tate $k$-Banach algebra over $\mathcal{A}$ of multiradius $\pmb{r}=(r_1,\dots, r_n)\in (\R_+)^N$ is the algebra over $k$
\[\Big\{\sum_{J\in \N^n}a_J\pmb{T}^J, \text{ }a_J\in A\text{ and }\lim_{|J|\to \infty}\algnorm{a_J}\cdot \pmb{r}^J=0 \Big\}\]
(for $J=(j_1,\dots, j_n)\in \N^n$, we denote $\prod_{i\in \card{n}}T_i^{j_i}$ by $\pmb{T}^J$ and $\prod_{i\in \card{n}}r_i^{j_i}$ by $\pmb{r}^J$) with a complete $k$-algebra norm defined by
\[\Big\vvvert\sum_{J\in \N^n}a_J\pmb{T}^J\Big\vvvert_{\mathcal{T}_{\mathcal{A}}(\pmb{r})}:=\sup_{J}\algnorm{a_J}\cdot \pmb{r}^J\]
This Banach algebra is denoted by $\mathcal{A}\{ r_1^{-1}T_1, \dots, r_n^{-1}T_n\}$, and is called an $\mathcal{A}$-Tate algebra of multiradius $\pmb{r}$.
\end{example}
\begin{definition}
    Let $\mathcal{A}_1, \mathcal{A}_2$ be two Banach $k$-algebras, and $\phi: A_1\to A_2$ be a homomorphism of $k$-algebras. We say that $\phi$ is a \emph{homomorphism of Banach $k$-algebras} if it is bounded as a $k$-linear map. A homomorphism of Banach $k$-algebra $\phi$ is often denoted by $\phi: \mathcal{A}_1\to \mathcal{A}_2$. A homomorphism of Banach $k$-algebra $\phi$ is called an \emph{isomorphism of Banach $k$-algebras} if there exists a homomorphism of Banach $k$-algebras $\psi:\mathcal{A}_2\to\mathcal{A}_1$ such that $\phi\circ \psi=\mathrm{Id}_{\mathcal{A}_2}$ and $\psi\circ \phi=\mathrm{Id}_{\mathcal{A}_1}$.
\end{definition}

\subsubsection{Spectrum}

Let $\mathcal{A}=(A, \algnorm{\ndot})$ be a Banach $k$-algebra. Let $\algnorm{\ndot}'$ be a $k$-algebra seminorm on $A$. One says that $\algnorm{\ndot}'$ is \emph{bounded} (with respect to $\mathcal{A}$) if there exists $C>0$ such that $\algnorm{\ndot}'\leq C\algnorm{\ndot}$. Its null-space is a closed ideal $I$ of $A$; the quotient $k$-algebra norm of $\algnorm{\ndot}'$ on the quotient $k$-algebra $A/I$ is bounded with respect to the quotient $k$-algebra norm of $\algnorm{\ndot}$. (\cite[Remark 1.2.2.i]{Ber})

\begin{definition}
Let $\mathcal{A}$ be a $k$-Banach algebra. The \emph{Berkovich spectrum} $\M(\mathcal{A})$ is the following topological space: the points, denoted by $z$, are bounded multiplicative $k$-algebra seminorms $\algnorm{\ndot}_z$ on $\mathcal{A}$, and the topology is the weakest topology on this set of points, for which all $\R_{\geq 0}$-valued functions of the form $z\mapsto \algnorm{f}_z$ are continuous for any $f\in A$. This topology is called the \emph{canonical topology}. For any subset $V$ of $\M(\mathcal{A})$, we denote by $\Int^{\mathrm{top}}(V)$ the topological interior of $V$. This topological interior is to be compared with the notion of \emph{interior of an affinoid subdomain in an affinoid domain} (see \cite[Definition 2.5.7]{Ber}), which we do not use in this article.
\end{definition}
\begin{remark}
    A basis for the canonical topology constituting of open sets is given by \emph{basic open sets}, which are sets of the form 
    \[U(f; p, q):=\{z\in \M(\mathcal{A}): p<\abs{f}_z<q \}\]
    indexed by $(p, q)\in \R^2$ and $f\in \mathcal{A}$. A general open set is a union of finite intersections of basic open sets.
\end{remark}

\begin{proposition}\label{banach spectrum compact Hausdorff}
    Let $\mathcal{A}$ be a $k$-Banach algebra. Then $\M(\mathcal{A})$ is a non-empty compact Hausdorff topological space. (\cite[Theorem 1.2.1]{Ber})
\end{proposition}

For any point $z\in \M(A)$, let $\mathfrak{p}_z$ be the closed ideal $\mathfrak n(\algnorm{\ndot}_z)$ of $A$ which is a prime ideal, and $f(z)$ be the image of $f$ in the quotient $k$-algebra $A/\mathfrak{p}_z$. The \emph{residual field} at $z$ is defined to be the fraction field of $A/\mathfrak{p}_z$, denoted by $\kappa(z)$, it is equipped with a quotient norm $\abs{\ndot}_z$ of $\algnorm{\ndot}_z$, which becomes an absolute value on $\kappa(x)$ extending $\abs{\ndot}$ on $k$. The \emph{completed residual field} at $z$ is defined to be the completion of $\abs{\ndot}_z$ with respect to this quotient norm $\abs{\ndot}_z$, denoted as $\widehat{\kappa}(z)$. The canonical homomorphism of $k$-algebra from $A$ to $(\widehat{\kappa}(z), \abs{\ndot}_z)$ is denoted by $\chi_z$. It is a homomorphism of Banach $k$-algebras.

\begin{definition}
Let $\mathcal{A}$ be a Banach $k$-algebra. A \emph{character} $\chi$ of $\mathcal{A}$ is a homomorphism of Banach $k$-algebra from $\mathcal{A}$ to some complete valued field extension $(K, \abs{\ndot}_K)$ of $(k, \abs{\ndot}_k)$. Two characters $\chi_1: \mathcal{A}\to (K_1,\abs{\ndot}_{K_1})$ and $\chi_2: \mathcal{A}\to (K_2,\abs{\ndot}_{K_2})$ are said to be \emph{equivalent} if there exist a character $\chi: \mathcal{A}\to (K,\abs{\ndot}_{K})$ and valued field extensions $\iota_1: K\to K_1$ and $\iota_2: K\to K_2$ which preserve norms such that $\chi=i_1\circ \chi_1=i_2\circ \chi_2$. Let $[\chi]$ be the equivalence class of $\chi$.
\end{definition}

\begin{lemma}
The set of points of $\M(\mathcal{A})$ is in canonical bijection with the set of equivalence classes of characters on $\mathcal{A}$. This bijection sends $z\in \M(\mathcal{A})$ to $[\chi_z]$. (\cite[Remark 1.2.2.ii]{Ber})
\end{lemma}

\begin{definition}
The \emph{Gelfand transform} of $\mathcal{A}$ is the homomorphism of Banach $k$-algebras
\[\widehat{}: A\to \prod_{z\in \M(\mathcal{A})}\hat{\kappa}(z),\quad f\mapsto\widehat{f}=(f(z))_{z\in \M(\mathcal{A})}\]
\end{definition}

\begin{proposition}\label{invertible spectrum}
An element $f\in \mathcal{A}$ is invertible if and only if $f(z)\neq 0$ for any $z\in \M(\mathcal{A})$. (\cite[Corollary 1.2.4]{Ber})
\end{proposition}

\subsubsection{Continuous map}

\begin{proposition}\label{morphism of spectrum}
    Let $\phi:\mathcal{A}_1\to \mathcal{A}_2$ be a homomorphism of Banach $k$-algebras. It induces a continuous map $\phi^{\star}: \M(\mathcal{A}_2)\to \M(\mathcal{A}_1)$ by sending an equivalent class of characters $[\chi]$ of $\mathcal{A}_2$ to the class of characters $[\chi\circ \psi]$ of $\mathcal{A}_1$. (\cite[Remark 1.2.2 (iii)]{Ber})
\end{proposition}

\begin{lemma}\label{injective}
    If $\phi:\mathcal{A}_1\to \mathcal{A}_2$ is a homomorphism of Banach $k$-algebras with dense image, then $\phi^{\star}$ is an injective map whose image is closed.
\end{lemma}
\begin{proof}
    The map $\phi^{\star}$ is injective since for any two characters $\chi_1, \chi_2: \mathcal{A}_2\to K$, if $\chi_1\circ \phi=\chi_2\circ \phi$, then the restriction of $\chi_1$ and $\chi_2$ on the image of $\phi$ are equal, hence the two characters are equal by the density of image.\par
    Let $(z, \abs{\ndot}_z)\in \M(\mathcal{A}_1)$ which is not in the image of $\phi^{\star}$, then $\ker(\phi)\nsubseteq \mathfrak{p}_z$: otherwise the character $\mathcal{A}_1/\ker(\phi)\to \hat{\kappa}(z)$ extends to a character $\mathcal{A}_2\to \hat{\kappa}(z)$ by the density of image of $\phi$. Now there exists $f\in \ker(\phi)\setminus\mathfrak{p}_z$, so $|f|_z\neq 0$. For small enough $\epsilon>0$, the basic open set $U(f;|f|_z-\epsilon,|f|_z+\epsilon)\subset \M(\mathcal{A}_1)$ is a neighbourhood of $(z, \abs{\ndot}_z)$ which is not contained in the image of $\phi^{\star}$. So the image of $\phi^{\star}$ is a closed subset in $\M(\mathcal{A}_1)$.
\end{proof}

\subsubsection{Spectral seminorm}

\begin{definition}
The \emph{spectral algebra seminorm} $\algnorm{\ndot}_{\mathrm{sp}}$ of an algebra seminorm $\algnorm{\ndot}$ on a $k$-algebra $A$ is the one defined by 
\[\forall f\in \mathcal{A}, \quad \algnorm{f}_{\mathrm{sp}}:=\lim_{\begin{subarray}{c}n\to \infty\end{subarray}}\algnorm{f^n}^{\frac{1}{n}}.\] 
Note that the triangle inequality for $\algnorm{\ndot}_{\mathrm{sp}}$ follows from sub-multiplicativity of $\algnorm{\ndot}$. In general, $\algnorm{\ndot}_{\mathrm{sp}}$ is only a seminorm even if $\algnorm{\ndot}$ is a norm.
\end{definition}
\begin{remark}
    The existence of limit is guaranteed by the (multiplicative) Fekete lemma for the sub-multiplicative sequence $\{\algnorm{f^n}\}_{n\in \N}$. The spectral seminorm is sub-multiplicative, and is bounded by the original seminorm $\algnorm{\ndot}$. Moreover, it is power-multiplicative by construction. 
\end{remark}

\begin{proposition}\label{Gelfand uniform}
    Let $\mathcal{A}$ be a $k$-Banach algebra. For any $f\in A$, one has  (\cite[Theorem 1.3.1]{Ber})
    \[\algnorm{f}_{\mathrm{sp}}=\max_{z\in \M(A)}|f|_z\]
\end{proposition}

\begin{definition}
Let $\mathcal{A}$ be a Banach $k$-algebra. The \emph{radical} of $\mathcal{A}$ is the null-space of its spectral seminorm $\mathfrak n(\algnorm{\ndot}_{\mathrm{sp}})$. A Banach $k$-algebra with radical equal to $\{0\}$ is called \emph{semi-simple}. Elements in the radical are said to be \emph{quasi-nilpotent} (or \emph{topological nilpotent}).
\end{definition}
\begin{remark}
    The radical of $\mathcal{A}$ contains the nil-radical of $A$; in other words, nilpotent elemtents are quasi-nilpotent. If $\mathcal{A}$ is semi-simple, then $A$ is reduced. The converse may not be true.
\end{remark}

Let $\mathcal{A}=(A, \algnorm{\ndot})$ be a $k$-Banach algebra. The spectral seminorm $\algnorm{\ndot}_{\mathrm{sp}}$ defines a quotient norm on the quotient $k$-algebra $\mathcal{A}/\rad(\mathcal{A})$, still denoted by $\algnorm{\ndot}_{\mathrm{sp}}$. The quotient norm $\algnorm{\ndot}_{\mathrm{sp}}$ is bounded by the quotient norm of $\algnorm{\ndot}$. The \emph{uniformization} $\mathcal{A}^u$ of $\mathcal{A}$ is defined to be the Banach $k$-algebra of separated completion of $(\mathcal{A}/\rad(\mathcal{A}), \algnorm{\ndot}_{\mathrm{sp}})$. Conversely, if $\algnorm{\ndot}$ is a power-multiplicative Banach algebra norm on $A$ with radical $\{0\}$, then it is said to be \emph{uniform}.

Obviously, $\algnorm{\ndot}_{\mathrm{sp}}$ is bounded by $\algnorm{\ndot}$. It is important to note that the converse may \textbf{not} be true in general. In other words, $\algnorm{\cdot}_{\mathrm{sp}}$ may \textbf{not} be complete on $\mathcal{A}/\text{rad}(\mathcal{A})$. Yet one still has the following statement
\begin{proposition}\label{uniformization keeps spectrum}
    $\M(\mathcal{A})$ is canonically homeomorphic to $\M(\mathcal{A}^u)$. (\cite[Corollary 1.3.3, 1.3.4]{Ber})
\end{proposition}

\subsubsection{Banach module}
One can also consider seminorms on modules over Banach algebra. Let $\mathcal{A}$ be a Banach $k$-algebra. A \emph{(semi)normed $\mathcal{A}$-module} is defined to be an $A$-module $M$ with a (semi)norm $\norm{\ndot}$ such that $(M, \norm{\ndot})$ is a (semi)normed vector space over $k$ (denoted by $\mathcal{M}$), and that the multiplication is bounded, in the sense that there exists $C>0$ such that
\[\forall a\in \mathcal{A},\  \forall m\in M, \quad \norm{a\cdot m}\leq C\algnorm{a}\cdot \norm{m}\]
One calls a \emph{Banach $\mathcal{A}$-module} a normed $\mathcal{A}$-module $(M, \norm{\ndot})$ whose norm is complete. 

Let $\mathcal{M}_1=(M_1, \norm{\ndot}_1)$, $\mathcal{M}_2=(M_2, \norm{\ndot}_2)$ be Banach $\mathcal{A}$-modules and $\phi: M_1\to M_2$ be a homomorphism of $A$-modules. It is called \emph{bounded} if there exists $C>0$ such that
$\norm{\phi(m_1)}_2\leq C\norm{m_1}_1$ for any $m_1\in M_1$. In this case $\phi$ is said to be a \emph{homomorphism of Banach $\mathcal{A}$-modules}, and is denoted by $\phi: \mathcal{M}_1\rightarrow \mathcal{M}_2$. In addition, the homomorphism $\phi$ of Banach $\mathcal{A}$-modules is called \emph{admissible} if it is admissible as linear map between normed-vector spaces over $k$.

\begin{definition}
    Let $\mathcal{M}$ be a Banach $\mathcal{A}$-module. It is called a \emph{Banach finite} $\mathcal{A}$-module if there exists $l\in \N_+$ and a surjective homomorphism of Banach $\mathcal{A}$-modules $\mathcal{A}^{\oplus l}\to\mathcal{M}$ where $\mathcal{A}^{\oplus l}$ is the Banach $\mathcal{A}$-module corresponding to the $A$-module $A^{\oplus l}$ equipped with the norm $(a_1, \dots, a_l)\mapsto \max\algnorm{a_i}$. (Note that such a homomorphism is necessarily admissible.)
\end{definition}

\begin{proposition}\label{Banach finite module}
    Let $\mathcal{A}$ be a Banach $k$-algebra and $\mathcal{M}$ be a Banach $\mathcal{A}$-module. If $A$ is Noetherian as a $k$-algebra and $M$ is finitely generated as $A$-module, then any $\mathcal{A}$-sub-module of $\mathcal{M}$ is closed, and $\mathcal{M}$ is a Banach finite $\mathcal{A}$-module. (\cite[Lemma 1.2.3]{FvdP} 
\end{proposition}

\begin{definition}
    Let $\phi:\mathcal{A}_1\to \mathcal{A}_2$ be a homomorphism between Banach $k$-algebras. It is called \emph{Banach finite} if $\mathcal{A}_2$ is a Banach finite $\mathcal{A}_1$-module. In this case $\mathcal{A}_2$ is called a Banach finite $\mathcal{A}_1$-algebra.
\end{definition}

\begin{remark}\label{Banach module finite surjective is finite}
    If a $k$-Banach algebra homomorphism $\phi$ is finite as homomorphism of $k$-algebra, and $\mathcal{A}_1$ is Noetherian, then $\phi$ is automatically Banach finite: there is a surjective $\mathcal{A}_1$-module homomorphism $p: \mathcal{A}_1^{\oplus n}\to \mathcal{A}_2$, by Proposition \ref{Banach finite module} $\ker(p)$ is closed. Then $p$ is continuous hence is admissible by Corollary \ref{equivalent norm}. So $\mathcal{A}_2$ is a Banach finite $\mathcal{A}_1$-module.
\end{remark}

\subsection{Affinoid algebras}\hfill\break

Affinoid algebras is a special kind of $k$-Banach algebras possessing good finiteness properties. These features allows one to endow a locally ringed space structure on their Berkovich spectra, namely the affinoid spaces. As a consequence, the Banach algebra norm of an affinoid algebra is equivalent to its spectral seminorm whenever the later is actually a norm. 

\subsubsection{Basic constructions}
Affinoid algebras are $k$-Banach algebras that are quotient algebras of Tate algebras. Among them are strict affinoid algebras which have good finiteness properties such as Noetherianity. Some good properties pass to general affinoid algebra by a technique enlarging the base valued field which makes the affinoid algebra strict.

\begin{definition}\label{standard Tate}
    For a multi-radius $\pmb{r}=(r_1,\dots, r_n)\in \R^n$, the algebra 
    $$k\{ r_1^{-1}T_1,\dots,r_n^{-1}T_n\}=\{f=\sum_{J\in \N^n}^{\infty}a_J\pmb{T}^{J}:a_{J}\in k, \abs{a_{J}}\pmb{r}^{J}\to 0 \text{ as } |J|\to \infty \}$$
    is called the \emph{Tate algebra} over $k$ with multi-radius $\pmb{r}$. Denote it by $\mathcal{T}_n(\pmb{r})$. It is a $k$-Banach algebra with respect to the \emph{Gauss norm of multi-radius $\pmb{r}$} defined by
    \[\algnorm{f}_{\mathcal{T}_n(\pmb{r})}=\max_{J}\abs{a_{J}}\pmb{r}^{J}\]
    One can define Tate algebra over other complete ultra-metric valued fields.
\end{definition}
\begin{remark}
    This Gauss norm is obviously sub-multiplicative. It is in fact multiplicative by an argument as in the proof of Gauss Lemma.
\end{remark}

\begin{definition}
    A $k$-Banach algebra $\mathcal{A}$ is called an \emph{affinoid algebra} if there exists an admissible surjective homomorphism from some Tate algebra $k\{\pmb{r}^{-1}\pmb{T}\}$ to $\mathcal{A}$. The Banach algebra norm on an affinoid algebra $\mathcal{A}$ is called an \emph{affinoid algebra norm}. If one can take $\pmb{r}$ with $r_i=1$ for all $i\in \card{n}$, then $\mathcal{A}$ is called a \emph{strict} affinoid algebra. One may define affinoid algebra similarly over other complete ultrametric valued field.
\end{definition}
\begin{remark}
    An affinoid algebra norm and the quotient algebra norm of Gauss algebra norm by the defining admissible surjective homomorphism are just equivalent but not necessarily equal.
\end{remark}

One can construct new affinoid algebras out of old ones by various algebraic operations.
\begin{example}\label{quotient of affinoid is affinoid}
    The quotient Banach algebra of an affinoid algebra is an affinoid algebra.
\end{example}

\begin{proposition}\label{finite extension is affinoid}
    Let $\mathcal{C}$ be a $k$-Banach algebra which is finite over an affinoid algebra $\mathcal{A}$, then $\mathcal{C}$ itself is an affinoid algebra. If $\mathcal{A}$ is strict, then $\mathcal{C}$ is strict.
\end{proposition}
\begin{proof}
    Let $\{c_i\}_{i\in \card{m}}\subset \mathcal{C}$ be a finite set of generators of $\mathcal{C}$ over $\mathcal{A}$, then consider an $\mathcal{A}$-Tate algebra $\mathcal{A}\{\pmb{r}^{-1} \pmb{T}\}$ where $r_i\geq \algnorm{c_i}_{\mathcal{C}}$. There is a surjective $k$-algebra homomorphism defined by
    \[\gamma: \mathcal{A}\{\pmb{r}^{-1} \pmb{T}\}\to \mathcal{C}, \text{ }T_i\mapsto c_i\]
    which is bounded as there exists $C>0$ such that
    $$\algnorm{\gamma(\sum_{J}a_J\pmb{T}^J)}_{\mathcal{C}}\leq \max_{J\in \N^{m}}\algnorm{a_J\pmb{c}^J}_{\mathcal{C}}\leq C \max_{J\in \N^m}\algnorm{a_J}_{\mathcal{A}}\cdot\pmb{r}^J$$
    By Corollary \ref{equivalent norm}, $\gamma$ is admissible, the norm $\algnorm{\ndot}_{\mathcal{C}}$ is equivalent to the quotient norm of the $\mathcal{A}$-Tate norm. Hence $\mathcal{C}$ is an affinoid algebra. The strictness is obtained by choosing $r_i\in \abs{k^{\times}}$ (see Lemma \ref{Tate strict}).
\end{proof}

\begin{proposition}\label{affinoid box}
    Let $\mathcal{B}$ be a Banach $k$-algebra. Suppose that there exists a finitely generated $k$-algebra $A$ which is dense in $\mathcal{B}$, then there exists an affinoid algebra $\mathcal{A}$ in which $A$ is a dense $k$-sub-algebra and a homomorphism of Banach $k$-algebras $\mathcal{A}\rightarrow \mathcal{B}$ which extends the identiy homomorphism on $A$.
\end{proposition}
\begin{proof}
    Let $\parenth{a_i}_{i\in \card{m}}$ be a set of generators of $A$. For each $i\in \card{m}$, let $r_i$ denote $\algnorm{a_i}_{\mathcal{B}}$ and let $\pmb{r}$ denote the multi-radius consisting of $\parenth{r_i}_{i\in \card{m}}$. Consider the Tate algebra $\mathcal{T}_{\pmb{r}}$ and the homomorphism of $k$-algebras
    \[k[T_1, \dots, T_m]\rightarrow \mathcal{B}, \quad T_i\mapsto a_i\]
    By the ultra-metricity of $\algnorm{\ndot}_{\mathcal{B}}$ and the definition of $\pmb{r}$, one has
    \[\forall
     n\in \N, \forall
     J\in \N^m, \forall f_J\in k, \algnorm{\sum_{J\in \N^m}f_J\cdot \pmb{a}^J}_{\mathcal{B}}\leq \algnorm{\sum_{J\in \N^m}f_J\cdot \pmb{a}^J}_{\mathcal{T}_{\pmb{r}}}\]
    so by a density argument one can extend it to a homomorphism of Banach $k$-algebras
    \[\mathcal{T}_{\pmb{r}}\rightarrow \mathcal{B}, \quad T_i\mapsto f_i\]
    Let $\mathscr{I}$ be the kernel ideal of this homomorphism. To conclude it suffices to take $\mathcal{A}$ as $\mathcal{T}_{\pmb{r}}/\mathscr{I}$.
\end{proof}

\begin{proposition}\label{Affinoid sub-algebra dense finite is finite}
    Let $(B, \algnorm{\ndot})$ be a normed algebra and let $\mathcal{B}$ be its separated completion. Let $A$ be a sub-$k$-algebra of $B$, equipped with the restriction algebra norm of $\algnorm{\ndot}$, and let $\mathcal{A}$ be the separated completion of $(A, \algnorm{\ndot})$. Assume that $\mathcal{A}$ is an affinoid algebra. If $B$ is integral and is finite over $A$, then $\mathcal{B}$ is Banach finite over $\mathcal{A}$. Therefore $\mathcal{B}$ is an affinoid algebra.
\end{proposition}
\begin{proof}
    By assumption, there exists $j\in \N$ and a homomorphism of $k$-algebras and elements $\parenth{e_i}_{i\in \card{j}}\subseteq B$ such that
    \[F: \bigoplus_{i\in \card{j}}A\rightarrow B, 1_i\mapsto e_i\]
    Moreover, $F$ is bounded
    \[\algnorm{\sum_{i\in \card{j}}a_i\cdot e_i}\leq \max_{i\in \card{j}}\algnorm{a_i\cdot e_i}\leq \max_{i\in \card{j}}\algnorm{e_i}\cdot \max_{i\in \card{j}}\algnorm{a_i}\]
    So $F$ extends to a homomorphism of Banach $\mathcal{A}$-modules
    \[\mathcal{F}: \bigoplus_{i\in \card{j}}\mathcal{A}\rightarrow \mathcal{B}, 1_i\mapsto e_i\]
    Let $\mathcal{B}^{-}$ be the image of $\mathcal{F}$, it is a Banach finite $\mathcal{A}$-module with the quotient norm $\norm{\ndot}_{\mathcal{F}}$ induced by $\mathcal{F}$. As $\mathcal{B}^{-}$ is Banach finite over $\mathcal{A}$, it is an affinoid algebra with an affinoid algebra spectral norm $\algnorm{\ndot}^{-}$, which is equivalent to $\norm{\ndot}_{\mathcal{F}}$. Now on $\mathcal{B}^{-}$, $\algnorm{\ndot}$ is bounded with respect to $\algnorm{\ndot}^{-}$ by the continuity of $\mathcal{F}$. To show the reverse, note that $\mathcal{B}^{-}$ is dense in $\mathcal{B}$, so by Theorem \ref{Gelfand uniform} one has for any $b\in \mathcal{B}^{-}$
    \[\algnorm{b}^{-}=\max_{z\in \M(\mathcal{B}^{-})}\abs{b(z)}=\max_{z\in \M(\mathcal{B})}\abs{b(z)}=\algnorm{b}_{\mathrm{sp}}\leq \algnorm{b}\]
    Therefore $\algnorm{\ndot}$ and $\algnorm{\ndot}^{-}$ are equivalent norms on $\mathcal{B}^{-}$, so $\mathcal{B}^{-}$ is closed in $\mathcal{B}$, hence coincides with it.
\end{proof}

To make an affinoid algebra strict, one can enlarge the base field.
\begin{lemma}\label{Tate field}
    Let $\pmb{r}=(r_1,\dots, r_n)$ be a multi-radius such $\{\alpha(\log r_i)\}_{i\in \card{n}}$ are $\Q$-linearly independent. Then the $k$-affinoid algebra
    $$K_{\pmb{r}}:=k\{\pmb{r}^{-1}\pmb{T}, \pmb{r}\pmb{T}^{-1}\}=k\{\pmb{r}^{-1}\pmb{T}, \pmb{r}\pmb{S}\}/(T_1S_1-1,\dots, T_nS_n-1)$$
    is a field. (\cite[Definition 2.1.1]{Ber})
\end{lemma}

\begin{lemma}\label{Tate strict}
    Let $\mathcal{T}_n(\pmb{r})$ be a $k$-Tate alegbra. It is strict if and only if $r_i\in \sqrt{\abs{k^{\times}}}$ for all $i$. (\cite[Corollary 2.1.6]{Ber})
\end{lemma}

\begin{corollary}\label{become strict Tate}
    Let $\mathcal{T}_n(\underline{r})$ be a $k$-Tate alegbra. Let $I\subseteq \card{n}$ be a subset of indices such that $\{\alpha(\log r_i)\}_{i\in I}$ are $\Q$-linearly independent and $|I|$ is maximal for this independence property. Let $\pmb{r}_I=(r_{i_1},\dots, r_{i_1})$, then $K_{\pmb{r}_I}\widehat{\otimes}_k \mathcal{T}_n(\pmb{r})$ is a strict $K_{\pmb{r}_I}$-Tate algebra.
\end{corollary}

\begin{corollary}\label{become strict affinoid}
    For any $k$-affinoid algebra $\mathcal{A}$, there exists a multi-radius $\pmb{r}_I=(r_i)_{i\in I}$ such that $\{\alpha(\log r_i)\}_{i\in I}$ are $\Q$-linearly independent and $K_{\pmb{r}_I}\widehat{\otimes}_k\mathcal{A}$ is a $K_{\pmb{r}_I}$-strict affinoid algebra. (\cite[Proposition 2.1.2]{Ber})
\end{corollary}

\subsubsection{Algebraic structures: Noetherianity}

Let $\mathcal{A}$ be a Banach $k$-algebra, one denotes by $\mathcal{A}^{\circ}$ the $k^{\circ}$-algebra $\{f\in\mathcal{A}\text{ }|\text{ }\algnorm{f}_{\mathcal{A},\text{sp}}\leq 1\}$, and by $\mathcal{A}^{\circ \circ}$ the ideal of $\mathcal{A}^{\circ}$ constituting of elements $\algnorm{f}_{\mathcal{A},\text{sp}}< 1$. The $\widetilde{k}$-algebra $\mathcal{A}^{\circ}/\mathcal{A}^{\circ\circ}$ is called the \emph{reduction} of $\mathcal{A}$. It can be shown that $\widetilde{\mathcal{T}_n}$ is isomorphic to $\widetilde{k}[T_1,\dots,T_n]$. (\cite[Proposition 5.1.2.2]{BGR})

\begin{definition}
    An element $f\in \mathcal{T}_n$ with $\algnorm{f}_{\mathcal{T}_n}=1$ is said to be \emph{regular in $z_n$ of degree $d$} if its reduction $\Bar{f}=\lambda (z_n)^d+\sum_{0\leq i\leq d-1}c_i(z_n)^{d-i}$ in $\Bar{\mathcal{T}_n}$ where $\lambda\in k^{\times}$ and $c_i\in \Bar{k}[z_1,\dots, z_{n-1}]$.
\end{definition}

\begin{proposition}\label{Weierstrass division}[Weierstrass division]
    Let $\mathcal{T}_n$ be the $k$-Tate algebra of multiradius $\underline{r}=\underline{1}$, then
    \begin{enumerate}[label=\rm{(\arabic*)}]
    \item Let $f\in \mathcal{T}_n$ be an distinguished element in $z_n$ of degree $d$, and $g\in \mathcal{T}_n$ be any element. Then there exist unique $r\in \mathcal{T}_{n-1}[z_n]$ of degree less than $d$ in $z_n$ and $q\in \mathcal{T}_n$ such that $g=q\cdot f+r$. Moreover $\algnorm{g}_{\mathcal{T}_n}=\max\{\algnorm{q}_{\mathcal{T}_n},\algnorm{r}_{\mathcal{T}_n}\}$
    \item Let $f\in \mathcal{T}_n$ with $\algnorm{f}_{\mathcal{T}_n}=1$. Then there exists a $k$-algebra automorphism $\tau$ of $\mathcal{T}_n$ such that $\tau(f)$ is regular in $z_n$.
    \end{enumerate}
    (\cite[Theorem 5.2.1.2]{BGR}, \cite[Theorem 3.1.1]{FvdP})
\end{proposition}

\begin{proposition}\label{Tate Noether closed ideal}
    The Tate algebra $\mathcal{T}_n$ is Noetherian. All of its ideals are closed. (\cite[Theorem 5.2.6.1, Corollary 5.2.7.2]{BGR}, \cite[Theorem 3.2.1]{FvdP})
\end{proposition}

\begin{corollary}\label{Affinoid Noether closed ideal}
    Any strict affinoid algebra is Noetherian. All of its ideals are closed (\cite[Proposition 6.1.1.3]{BGR}, \cite[Theorem 3.2.1]{FvdP}).
    Any affinoid algebra is Noetherian. All of its ideals are closed (\cite[Propositon 2.1.3]{Ber}).
\end{corollary}



\begin{proposition}\label{Noether normalization strict affinoid}[Noether normalization]
    For strict affinoid algebra $\mathcal{A}$, there exists an injective finite and admissible Banach algebra homomorphism $\mathcal{T}_d\to \mathcal{A}$ for some $d>0$. Moreover, $d$ equals the Krull dimension of $\mathcal{A}$. (\cite[Theorem 6.1.2.1]{BGR}, \cite[Theorem 3.2.1]{FvdP})
\end{proposition}

\begin{corollary}\label{Nulstellensatz}
    Let $\mathfrak{m}$ be a maximal ideal of strict affinoid algebra $\mathcal{A}$, then $\mathcal{A}/\mathfrak{m}$ is a finite extension of $k$.
\end{corollary}

\subsubsection{Topological structures: the spectral norm}

The Gauss norm on Tate algebra is equal to its spectral norm. For a general strict redueced affinoid algebra, the Banach algebra norm is equivalent to its spectral seminorm, thanks to the compatibility of Banach algebra norms with algebraic structures. 



One studies the spectral norm of the Tate algebra case by direct calculation.
\begin{proposition}\label{maximum principle Tate}
    For any $f\in \mathcal{T}_n$, there exists $z\in \mathrm{Max}(\mathcal{T}_n)$ such that $\abs{f(z)}_z=\algnorm{f}_{\mathcal{T}_n}$ (\cite[Proposition 5.1.4.3]{BGR}).
    On $\mathcal{T}_n$, the three norms are equal: $\algnorm{\ndot}_{\mathcal{T}_n}=\algnorm{\ndot}_{\mathcal{T}_n, \text{sp}}=\algnorm{\ndot}_{\mathcal{T}_n, \mathrm{spM}}$.
\end{proposition}

One then uses Noether normalization to investigate the spectral seminorm of general affinoid algebra.


\begin{proposition}\label{affinoid spectral norm complete}
    Let $\mathcal{A}$ be a reduced strict affinoid algebra. Then its spectral norm $\algnorm{\ndot}_{\mathcal{A},\text{sp}}$ is a complete norm on $\mathcal{A}$. It is equivalent to the Banach algebra norm $\algnorm{\ndot}_{\mathcal{A}}$. (\cite[Theorem 3.4.9]{FvdP}, \cite[Theorem 6.2.4.1]{BGR})
\end{proposition}

\begin{corollary}\label{affinoid norm equivalent to spectral norm}
    Let $\mathcal{A}$ be a reduced general affinoid algebra. Then there exists $C>0$ such that $\algnorm{f}\leq C\algnorm{f}_{\mathrm{sp}}$ for all $f\in \mathcal{A}$. In particular, $\algnorm{\ndot}_{\mathrm{sp}}$ is complete on $\mathcal{A}$ , and is equivalent to $\algnorm{\ndot}$. (\cite[Proposition 2.1.4.ii]{Ber})
\end{corollary}

\begin{remark}
    The constant $C$ here does not depend on $f\in\mathcal{A}$, it is uniform.
\end{remark}

\subsubsection{Affinoid space as locally ringed space}

The Berkovich spectrum of affinoid algebras are called affinoid spaces. It is possible to put locally ringed space structures on them. The construction of structural sheaf goes first with a Grothendieck topology generated by closed compact subsets of affinoid domains, then passes to the canonical topology by a limit process approximating an open set by these compact sets.

\paragraph{Affinoid domains and structural algebra}

\begin{definition}
    Let $\mathcal{A}$ be an affinoid algebra. An \emph{affinoid domain} is a closed subset $V$ of $\M(\mathcal{A})$, which is homeomorphic to $(\iota_V)^{\star}(\M(\mathcal{A}_V))$ for some affinoid algebra $\mathcal{A}_V$ and Banach algebra homomorphism $\iota_V: \mathcal{A}\to \mathcal{A}_V$, and satisfies the \emph{universal mapping property}: for any Banach algebra homomorphism $\phi: \mathcal{A}\to\mathcal{C}$ between affinoid algebras with $\phi^{\star}(\M(\mathcal{C}))\subseteq V$, there exists a unique Banach algebra homomorphism $\psi: \mathcal{A}_V\to \mathcal{C}$ with $\phi=\psi\circ\iota_V$
\end{definition}

\begin{lemma}
    Let $V$ be an affinoid domain in $\M(\mathcal{A})$. Then $V$ is homeomorphic to $\M(\mathcal{A}_V)$. Moreover $\mathcal{A}_V$ is a flat $\mathcal{A}$-algebra. (\cite[Proposition 2.2.4]{Ber})
\end{lemma}

\begin{example}
    Given $f=(f_1, \dots, f_m)$ and $g=(g_1, \dots, g_n)$ tuples of elements of $\mathcal{A}$, $p=(p_1,\dots, p_m)\in (\R_+^*)^m$ and $q=(q_1,\dots, q_n)\in (\R_+^*)^n$, the closed subset 
    $$V=\M(\mathcal{A})(p^{-1}f, qg^{-1}):=\{z\in \M(\mathcal{A}), |f_i(z)|_z\leq p_i, \text{ }|g_j(z)|_z\geq q_j\}$$
    is an affinoid domain. The corresponding homomorphism of affinoid algebras is
    $$\mathcal{A}\to \mathcal{A}_V=\mathcal{A}\{p_1^{-1}T_1, \dots, p_m^{-1}T_m, q_1S_1, \dots ,q_nS_n\}/(T_i-f_i, g_jS_j-1)$$
    Such domains are called \emph{Laurent domains}. If $n=0$, they are called \emph{Weierstrass domains}. 
\end{example}


\begin{lemma}\label{affinoid intersection}
    A finite intersection of affinoid domains is an affinoid domain. (\cite[Remark 2.2.2.iv]{Ber})
\end{lemma}

\begin{corollary}\label{affinoid topological basis}
    Any point $z\in \M(\mathcal{A})$ has a fundamental system of (closed) neighbourhoods consisting of affinoid domains. (\cite[Proposition 2.2.3]{Ber})
\end{corollary}

\paragraph{Special domains and acyclicity of structural presheaf}

\begin{definition}
    A \emph{special domain} $V$ in $\M(\mathcal{A})$ is a finite union of affinoid domains $V_i$ in $\M(\mathcal{A})$.
\end{definition}

\begin{definition}
    The \emph{Grothendieck topology} on $\M(\mathcal{A})$ is the one with special domains as admissible open sets and finite covering as admissible coverings. One notes $\M(\mathcal{A})_G$ for the space with this G-topology.
\end{definition}

\begin{definition}
    Let $\mathfrak{V}$ be an admissible covering of $\M(\mathcal{A})$ by affinoid domains $\parenth{V_i}_{i\in I}$, where $I$ is a finite set. Then for a Banach finite $\mathcal{A}$-module $\mathcal{M}$, the \emph{Cech complex} of $\mathcal{M}$ with respect to $V_i$ is defined to be the complex of Banach $\mathcal{A}$-modules
    $$C^{\centerdot}(\mathcal{M}, \mathfrak{V}):\text{ }0\to \mathcal{M} \to \prod_{i\in I}\mathcal{M}_i\to \prod_{i, j\in I}\mathcal{M}_{i, j}\to \dots$$
\end{definition}

One would like to have acyclicity of the complex $C^{\centerdot}(\mathcal{M}, \mathfrak{V})$ in order to follow standard construction of a structural sheaf on $\M(\mathcal{A})_{G}$.

\begin{theorem}\label{acyclicity strict affinoid}
    Let $\mathcal{A}$ be a strict affinoid algebra and $\mathfrak{V}$ an admissible covering by strict affinoid domains for $\M(\mathcal{A})$. Then $C^{\centerdot}(\mathcal{A}, \mathfrak{V})$ is acyclic. (\cite[Proposition 8.2.2.5]{BGR})
\end{theorem}

\begin{corollary}\label{acyclicity general affinoid}
    For general affinoid domain $\M(\mathcal{A})$ with general affinoid domains covering $\mathfrak{V}$, the complex $C^{\centerdot}(\mathcal{A}, \mathfrak{V})$ is acyclic. So is $C^{\centerdot}(M, \mathfrak{V})$ for finite Banach $\mathcal{A}$-module $M$. (\cite[Proposition 2.2.5]{Ber})
\end{corollary}

\begin{definition}\label{special domain}
    Let $V$ be any special domain in $\M(\mathcal{A})$. Fix a way of writing $V$ as $\bigcup_{i\in I}V_i$ where $I$ is a finite set and $V_i=\M(\mathcal{A}_{V_i})$ are affinoid algebras, let
    \[\mathcal{A}_V:=\ker(\prod_{i\in I}\mathcal{A}_{V_i}\to \prod_{i, j\in I}\mathcal{A}_{V_{i}\cap V_j})\]
    be the $k$-Banach algebra with sub-norm. The \emph{structural pre-sheaf of affinoid algebras} $\mathscr{O}_{\M(\mathcal{A})_G}$ on $\M(\mathcal{A})_G$ (\textbf{with respect to the G-topology}) is the one assigning $V$ the $k$-Banach algebra $\mathcal{A}_V$. It is a sheaf thanks to Corollary \ref{acyclicity general affinoid}.
\end{definition}

\begin{remark}
    The $k$-Banach algebra $\mathscr{O}_{\M(\mathcal{A})_G}(V)$ does not depend on the way of being a union of affinoid domains.
\end{remark}

\begin{definition}
    For any open subset $U$ of $\M(\mathcal{A})$, let $\mathscr{O}_{\M(\mathcal{A})}$ be the pre-sheaf of $k$-algebras (\textbf{with respect to the canonical topology}) which assigns $U$ the limit
    \[\mathscr{O}_{\M(\mathcal{A})}(U):=\varprojlim_{V\subset U, \text{V special domain}}\mathcal{A}_V\]
    It is also a sheaf thanks to the compactness of special domains under canonical topology. This is called the \emph{structural sheaf} of $\M(\mathcal{A})$.
\end{definition}

\begin{proposition}
    $\mathscr{O}_{\M(\mathcal{A})}$ is a sheaf of local rings. The topological space $\M(\mathcal{A})$ has a structure of locally ringed space given by the sheaf $\mathscr{O}_{\M(\mathcal{A})}$. (\cite[Section 2.3]{Ber})
\end{proposition}

\subsection{Spectral calculus}\hfill\break

Gelfand-Shilov theory allows one to do multi-variable spectral calculus for (commutative) Banach algebras over $\C$. In particular, one can localize a homomorphism between Banach algebras onto a neighbourhood of its spectrum. Similar theory, as develloped in \cite[Chapter 7]{Ber}, exists in the non-Archimedean base field setting.

\subsubsection{Holomorphic envelop}

The holomorphic convexity of spectrum of a homomorphism of Banach $k$-algebra depends on the dense-ness of its image. In case where the spectrum of a homomorphism is not holomorphic convex, one can add variables to the source algebra so that spectrum of extended homomorphism is holomorphically convex. 

\begin{definition}
Let $\mathcal{A}$ and $\mathcal{B}$ be Banach $k$-algebras, and $\phi: \mathcal{A}\to \mathcal{B}$ be a homomorphism of Banach algebras. The \emph{spectrum of homomorphism $\phi$} is the image of $\M(\mathcal{B})$ in $\M(\mathcal{A})$ under $\phi^{\star}$. Denote it by $\Sigma_{\phi}$
\end{definition}

\begin{definition}
Let $\mathcal{A}$ be a Banach $k$-algebra. Let $\Omega$ be a compact subset of $\M(\mathcal{A})$. The \emph{holomorphic convex envelop} of $\Omega$ in $\M(\mathcal{A})$ is the subset
\[\Omega^{\mathrm{h}}:=\parenth{z\in \M(\mathcal{A}) \ |\ \forall f\in \mathcal{A}, \abs{f}_z\leq \sup_{z'\in \Omega}\abs{f}_{z'}}\]
The subset $\Omega$ is said to be \emph{holomorphically convex} if $\Omega^{\mathrm{h}}=\Omega$.
\end{definition}

\begin{lemma}\label{Hull by Weierstrass}
    The intersection of all Weierstrass neighbourhoods of $\Omega$ in $\M(\mathcal{A})$ coincide with $\Omega^{\mathrm{h}}$. (\cite[Proposition 2.6.1]{Ber})
\end{lemma}

\begin{proposition}\label{holomorphic convex spectrum}
    Let $\mathcal{A}$ be a $k$-affinoid algebra, $\mathcal{B}$ be a Banach $k$-algebra. Let $\phi: \mathcal{A}\to \mathcal{B}$ be a homomorphism of Banach $k$-algebras. Let $\mathcal{B}'$ be the closed sub-algebra generated by the image of $\phi$ of $\mathcal{A}$ in $\mathcal{B}$ and let $\phi':\mathcal{A}\to \mathcal{B}'$ be the restricted homomorphism. Then $(\Sigma_{\phi})^{\mathrm{h}}=\Sigma_{\phi'}$. (\cite[Proposition 7.3.1]{Ber})
\end{proposition}

\begin{corollary}
    Let $\phi: \mathcal{A}\to \mathcal{B}$ be a homomorphism of Banach algebras from an affinoid algebra to a Banach algebra with dense image. Then $\Sigma_{\phi}$ is holomorphically convex.
\end{corollary}

One has the following analogue of Arens-Calderon theorem, which holomorphically convexifies the spectrum of a homomorphism of Banach $k$-algebras by adding variables on the source algebra.
\begin{proposition}\label{Arens-Calderon}
    Let $\mathcal{A}$ be a $k$-affinoid algebra, $\mathcal{B}$ be a Banach $k$-algebra. Let $\phi: \mathcal{A}\to \mathcal{B}$ be a homomorphism of Banach $k$-algebras. Then for any open neighbourhood $U$ in $\M(\mathcal{A})$ of the spectrum $\Sigma_{\phi}$, there exists a homomorphism of Banach algebras extending $\phi$
    $$\widetilde{\phi}: \widetilde{\mathcal{A}}:=\mathcal{A}\{r_1^{-1}T_1, \dots, r_n^{-1}T_n\} \to \mathcal {B}$$
    such that $\mathrm{pr}((\Sigma_{\phi})^{\mathrm{h}})\subseteq U$, where $\mathrm{pr}: \M(\widetilde{\mathcal{A}})\to \M(\mathcal{A})$ is the canonical map of projection. (\cite[Proposition 7.3.3]{Ber})
\end{proposition}

\subsubsection{Holomorphic functional calculus}

It is easy to localize the homomorphism to holomorphic convex neighbourhood of its spectrum. For a spectrum of homomorphism which is not holomorphically convex, one uses Proposition \ref{Arens-Calderon} to localize the homomorphism to any neighbourhood of it.

\begin{lemma}\label{polyhedra spectral localisation}
    Let $\phi:\mathcal{A}\to \mathcal{B}$ be a Banach algebra homomorphism from an affinoid algebra $\mathcal{A}$ to a Banach algebra $\mathcal{B}$. Then for any Laurent domain neighbourhood $V$ of $\Sigma_{\phi}$, $\phi$ extends to a unique Banach algebra homomorphism $\phi_V: \mathcal{A}_V\to \mathcal{B}$. (\cite[Corollary 2.5.16]{Ber})
\end{lemma}

\begin{theorem}\label{holomorphic functional calculus}
    Let $\phi: \mathcal{A}\to \mathcal{B}$ be a homomorphism of Banach algebras from an affinoid algebra to a Banach algebra. Let $V\subseteq\M(\mathcal{A})$ be any special domain containing $\Sigma_{\phi}$. Then there exists a Banach algebra homomorphism 
    $$\theta_{\phi}: \Gamma(V, \mathscr{O}_{\M(\mathcal{A})})\to \mathcal{B}$$ 
    satisfying $\phi=\theta_{\phi}\circ\iota_V$, where $\iota_V: \mathcal{A}\to \mathcal{A}_V=\Gamma(V, \mathscr{O}_{\M(\mathcal{A})})$ is the Banach algebra homomorphism corresponding to the inclusion $V\subseteq \M(\mathcal{A})$. (\cite[Theorem 7.3.4]{Ber})
\end{theorem}

\begin{remark}
    One can verify that the resulting Banach algebra homomorphism does not depend on the choice of $\widetilde{\phi}$.
\end{remark}

\subsection{Analytification of scheme of finite type}\hfill\break

There is a construction of Berkovich spectrum for a $k$-algebra similar to the one for $k$-Banach algebra, giving rise to analytification of $k$-schemes of locally finite type, as developped in \cite[Section 3.4]{Ber}.

\subsubsection{Local situation}
For affine varieties, the topological space of its analytification is defined in the same way as the spectrum of Banach algebra, except that boundedness requirement of seminorms are dropped. They enjoy similar basic properties as the spectrum of Banach algebra. Proofs are of same spirit hence are omitted.
\begin{definition}
    Let $Z=\spec (A_Z)$ be an affine $k$-scheme of finite type, where $A_Z$ is a $k$-algebra of finite type. Its Berkovich analytification $Z^{an}$ is the topological space constituting of all \emph{multiplicative} seminorms $\algnorm{\ndot}$ on $A_Z$ as points and with the canonical topology (the weakest topology making every function $\algnorm{\ndot}\to \algnorm{f}$ continuous for each $f\in A_Z$). (\cite[Remark 3.4.2]{Ber})
\end{definition}

\begin{definition}
    A \emph{character} on $A_Z$ is a homomorphism of $k$-algebra from $A_Z$ to some valued field extension $(K, \abs{\ndot}_K)$ over $(k, \abs{\ndot}_k)$. Two characters $\chi_1: A_Z\to K_1$ and $\chi_1: A_Z\to K_1$ are called \emph{equivalent} if there exists a $k$-algebra homomorphism $\chi_3: A_Z\to K_3$ and norm preserving $k$-algebra homomorphisms $i_1: K_1\to K_3$ and $i_2: K_2\to K_3$ satisfying $\chi_3=i_1\circ \chi_1=i_2\circ \chi_2$.
\end{definition}

\begin{lemma}
    There is a bijective map from the set of points of $Z^{an}$ to the set of equivalent classes of characters on $A_Z$.
\end{lemma}

\begin{proposition}\label{affine morphism analytification}
    Let $\phi: A_Z\to A_W$ be a homomorphism of $k$-algebras of finite type where $Z=\spec(A_Z)$ and $W=\spec(A_W)$. Then there is an induced continuous map $\phi^{\star}: W^{an}\to Z^{an}$, which sends a multiplicative seminorm $|\cdot|_w$ to $|\phi(\cdot)|_w$.
\end{proposition}

\begin{proposition}\label{injective and surjective morphism of non-normed spectrum}
    If $\phi$ is surjective, then $\phi^{\star}$ is injective and is a closed map; if $\phi$ is finite, then $\phi^{\star}$ is surjective. ( \cite[Proposition 3.46 (6)(7)]{Ber})
\end{proposition}

\begin{proposition}\label{normed spectrum into non-normed spectrum}
    Let $Z=\spec (A_Z)$ be an affine $k$-variety, $\algnorm{\ndot}$ an algebra norm on $A_Z$ and $\mathcal{A}_Z$ be the $k$-Banach algebra obtained by completing $A_Z$ with respect to $\algnorm{\ndot}$. Then the canonical homomorphism of $k$-algebras from $A_Z$ to $\mathcal{A}_Z$ induces a continuous map which embeds the Berkovich spectrum $\M(\mathcal{A}_Z)$ into $Z^{an}$ as a compact subspace (and is closed since $Z^{\mathrm{an}}$ is Hausdorff), and the Berkovich topology coincides with the induced topology from $Z^{an}$.
\end{proposition}
\begin{proof}
    For any $z\in \M(\mathcal{A}_Z)$, the multiplicative algebra seminorm (or the corresponding character) $\chi_z$ on $\mathcal{A}_Z$ corresponds to a unique multiplicative algebra seminorm on $A_Z$ by restriction. Since $A_Z$ is dense in $\mathcal{A}_Z$, the family of open sets $\{U(f; p, q), \text{ }f\in A_Z, \text{ } p,q \in \R\}$ form a basis for topology on $\M(\mathcal{A}_Z)$, hence the inherited topology coincides with the originial topology. So the embedding is continuous, and the image of $\M(\mathcal{A}_Z)$ is compact in $Z^{an}$. Since the topology on $Z^{an}$ is Hausdorff, the image of $\M(\mathcal{A}_Z)$ is closed.
\end{proof}

\begin{definition}
    An analytic function on open set $U\subseteq (\spec A_Z)^{an}$ is a map $h: U\to \coprod_{z\in U}\hat{\kappa}(z)$ which is a local uniform limit of rational functions: every $z\in U$ has an open neighbourhood $U'\subseteq U$ such that for every $\epsilon>0$, there exists $f_{U'}, g_{U'}\in A_Z$ with $|h(z)-\frac{f_{U'}(z)}{g_{U'}(z)}|<\epsilon$ and $g(z)\neq 0$ for all $z\in U'$. Denote by $\mathcal{R}^{\mathrm{an}}(U)$ the $k$-algebra of all analytic functions on $U$.
\end{definition}

\begin{definition}
    The structural sheaf $\mathscr{O}_{Z^{\mathrm{an}}}$ on $Z^{\mathrm{an}}$ is the one assigning $\mathcal{R}^{\mathrm{an}}(U)$ to an open set $U$. 
\end{definition}

\begin{proposition}
    $\mathscr{O}_{Z^{\mathrm{an}}}$ is a sheaf of local rings. The pair $(Z^{\mathrm{an}}, \mathscr{O}_{Z^{\mathrm{an}}})$ gives rise to a locally ringed space.
\end{proposition}

\begin{proposition}\label{morphism non-normed spectrum as ringed space}
    If $\mathcal{A}_Z$ is an affinoid algebra $\mathcal{A}_Z$, then there is a morphism of locally ringed space
    $$(\M(\mathcal{A}_Z, \mathscr{O}_{\M(\mathcal{A}_Z)})\to (Z^{an}, \mathscr{O}_{Z^{an}})$$
\end{proposition}
\begin{proof}
    The map of topological spaces is given in Proposition \ref{normed spectrum into non-normed spectrum}. For the ring homomorphism, it suffices to construct a $k$-algebra homomorphism $\mathscr{O}_{Z^{an}}(U)\to \mathcal{A}_V$ for any open set $U\subseteq \M(\mathcal{A})$ and any affinoid domain $V\subseteq U$. Moreover, it suffices to consider $U$ and $V$ of basic form
    \[U=U(\underline{p}^{-1}\underline{f}, \underline{q}\underline{g}^{-1}),\quad V=\M(\mathcal{A}_Z((\underline{p}-\underline{\epsilon})^{-1}\underline{f}, (\underline{q}+\underline{\epsilon})\underline{g}^{-1})) \text{ , }\epsilon>0\]
    There is a homomorphism of $k$-algebras $\mathcal{R}^{\mathrm{an}}(U)\to \mathcal{A}_Z((\underline{p}-\underline{\epsilon})^{-1}\underline{f}, (\underline{q}+\underline{\epsilon})\underline{g}^{-1})$ sending $\frac{\tilde{f}}{\tilde{g}}$ for $\tilde{f}, \tilde{g}\in A_Z$ to itself, the later being an element of $\mathcal{A}_V$ since $\frac{1}{\tilde{g}}\in \mathcal{A}_V$ by Lemma \ref{invertible spectrum}. As uniform limits of sequence in $\mathcal{R}^{\mathrm{an}}(U)$ remains to be uniform limits, this homomorphism extends to a $k$-algebra homomorphism $\mathscr{O}_{Z^{an}}(U)\to \mathcal{A}_V$.
\end{proof}

\subsubsection{Global situation}
One can analytify a scheme of finite type defined over $k$ by glueing local constructions.
\begin{definition}
    Let $X$ be a finite type scheme over $\spec k$, and write $X$ as $\bigcup X_i$ where $X_i=\spec A_{X_i}$ are affine charts. The Berkovich analytification of $(X, \mathscr{O}_X)$ is the locally ringed space obtained by gluing the Berkovich analytification $((X_i)^{an}, \mathscr{O}_{(X_i)^{an}})$ of each $(X_i, \mathscr{O}_{X_i})$. 
\end{definition}

\begin{proposition}\label{morphism analytification}
    Let $\phi: X\rightarrow Y$ be a morphism of schemes of locally finite type over $\spec k$. Then it induces a continuous map $\phi^{\mathrm{an}}: X^{\mathrm{an}}\rightarrow Y^{\mathrm{an}}$. And $\phi$ is (1) separated, (2) injective, (3) surjective, (4) an open immersion and (5) an isomorphism if and only if $\phi^{\mathrm{an}}$ has the same property. (\cite[Proposition 3.4.6]{Ber})
\end{proposition}

\begin{theorem}\label{Analytification Hausdorff compact}
    If $X$ is proper, then $X^{\mathrm{an}}$ is Hausdorff and compact. (\cite[Theorem 3.4.8]{Ber})
\end{theorem}


\section{Normed section algebra}

In this section, one studies norms on graded linear series of a line bundle on a projective variety.

\subsection{Basic setting}\hfill\break
Let $k$ be a field equipped with a complete non-Archimedean absolute value $\abs{\ndot}$, which is not trivial. Let $X$ be an irreducible scheme of finite type over $\spec k$. One denotes by $X^{\mathrm{an}}$ the Berkovich analytic space associated with $X$ and by $j_X: X^{\mathrm{an}}\rightarrow X$ the map sending any $x\in X^{\mathrm{an}}$ to its associated scheme point. 

\begin{enumerate}[fullwidth,itemindent=0em,label={\bf\arabic*.}]

\item\label{Item: algebra norm} Let $V=\bigoplus_{n\in \N}V_n$ be a graded $k$-algebra. Let $\algnorm{\ndot}$ be an algebra seminorm on $V$. For every $n\in \N$, this algebra seminorm induces by restriction a seminorm on the $k$-vector space $V_n$, denoted by $\norm{\ndot}_n$. As $\algnorm{\ndot}$ is sub-multiplicative, these seminorms satisfy the property
\[\forall (m,n)\in \N^2, s_m\in V_m, s_n\in V_n, \quad \norm{s_m\cdot s_n}_{m+n}\leqslant \norm{s_m}_m\cdot \norm{s_n}_n; \quad \norm{1}_0=1\]
Conversly, given a familly of \emph{ultrametric} seminorms $\{\norm{\ndot}_n\}_{n\in \N}$ on $k$-vector spaces $V_n$ of $V$ satisfying these properties, the seminorm on the graded $k$-algebra $\bigoplus_{n\in \N}V_n $ defined by 
\[\forall \underline{s}=(s_n)_{n\in \N},\quad \algnorm{\underline{s}}:=\sup_{n\in \N}\norm{s_n}_n\]
is submultiplicative, hence is an algebra seminorm on $V$. In fact, let $\underline{s}=(s_n)_{n\in\N}$ and $\underline{t}=(t_n)_{n\in\N}$ be two elements of $\bigoplus_{n\in \N}V_n$ and $\underline{u}=(u_n)_{n\in\N}=\underline{s}\cdot \underline{t}$, then one has
\[u_n=\sum_{\begin{subarray}{c}(p,q)\in\mathbb N^2\\
p+q=n\end{subarray}}s_p\cdot t_q.\]
By using the fact that the seminorm $\norm{\ndot}_{n}$ is ultrametric, one obtains that
\[\norm{u_n}_{n}\leqslant\max_{\begin{subarray}{c}
(p,q)\in\mathbb N^2\\
p+q=n
\end{subarray}}\norm{s_p\cdot t_q}_n\leqslant\max_{\begin{subarray}{c}
(p,q)\in\mathbb N^2\\
p+q=n
\end{subarray}}\norm{s_p}_{p}\cdot\norm{t_q}_{q},\]
so $\algnorm{\underline{u}}$ is bounded from above by $\algnorm{\underline{s}}\cdot\algnorm{\underline{t}}$. Denote by $\widehat{V}(\algnorm{\ndot})$ the separated completion of the seminormed algebra $(\bigoplus_{n\in \N}V_n, \algnorm{\ndot})$. 

One denotes by $\Upsilon(V_{\sbullet}(L))$ the set of all power-multiplicative ultrametric algebra norms $\algnorm{\ndot}$ which satisfies
\[\forall (s_n)\in V_{\sbullet}(L),\ \algnorm{(s_n)}=\sup_{n\in\N}\algnorm{s_n}.\]
This last condition is equivalent to the orthogonality of $\parenth{V_n}_{n\in \N}$ as $k$-linear subspaces.

\item\label{Item: algebra of sections} 
For any invertible $\mathscr{O}_X$-module $L$, one denotes by $V_\sbullet(L)$ the graded $k$-algebra $\bigoplus_{n\in\mathbb N}V_n(L)$ where $V_n(L):=H^0(X,L^{\otimes n})$. As $X$ is irreducible, $V_0(L)=k$. 


Let $f:Y\rightarrow X$ be a morphism of $k$-schemes. The morphism of $\mathscr{O}_X$-modules $L^{\otimes n}\to f_*f^*(L^{\otimes n})$ induces linear maps of $k$-vector spaces $V_n(L)\rightarrow V_n(f^*L)$ and graded homomorphism of degree $0$ of graded-$k$-algebras $V_\sbullet(L)\rightarrow V_\sbullet(f^*L)$. Denote by $V_n(L_{X|Y})$ and $V_{\sbullet}(L_{X|Y})$ the image vector space and image graded algebra.

One denotes by $Tot(L^{\vee})$ the scheme $\spec (\sym_{\mathscr{O}_X} L)$ over $\spec k$, by $\pi_X$ the canonical morphism of schemes $Tot(L^{\vee}) \rightarrow X$, and by $\mathbb{O}$ the reduced closed subscheme of zero section. If the graded $k$-algebra $V_{\sbullet}(L)$ is of finite type, then one denotes by $\pmb{0}$ the closed point of $\spec V_{\sbullet}(L)$ given by the maximal ideal $V_{\geq 1}(L)$, and by $p_X(\mathbf{0})$ the canonical blow-up morphism $Tot(L^{\vee})\rightarrow \spec V_{\sbullet}(L)$ along the sub-scheme $\mathbf{0}$.

Denote the integer $\dim_{k}V_n(L)-1$ by $d_n$. If $L^{\otimes n}$ is globally generated, there is a morphism induced by $V_n(L)$
\[\iota_n: \ X\rightarrow \mathbb{P}(V_n(L))\simeq\mathbb{P}^{d_n}.\]

\item\label{Item: metric} Let $L$ be an invertible $\mathscr{O}_X$-module. Let $\mathcal F_X$ be the sheaf of real-valued functions on $X^{\mathrm{an}}$. By \emph{pseudometric} on $L$ one refers to a morphism of sheaves of sets $\phi:L\rightarrow j_{X,*}(\mathcal F_X)$ such that, for any $x\in X^{\mathrm{an}}$, the map $\abs{\ndot}_\phi(x):L(x)\rightarrow\mathbb R$ induced by $\phi$ is a seminorm on the one-dimensional vector space $L(x)$ over $\widehat{\kappa}(x)$. If, for any $x\in X^{\mathrm{an}}$, the map $\abs{\ndot}_\phi(x)$ is a norm, one says that $\phi$ is a \emph{metric}.

We say that a pseudometric $\phi$ is \emph{(upper semi-)continuous} if, for any Zariski open subset $U$ of $X$ and any section $s\in\Gamma(X,L)$, the function $(x\in U^{\mathrm{an}})\rightarrow |s|_\phi(x)$ is (upper semi-)continuous. 


\item\label{Item: disc bundle}The pair $(L, \phi)$ is called a \emph{pseudometrized invertible $\mathscr{O}_X$-module}. For any $\epsilon\in\R$, the following subset of $Tot(L^{\vee})^{\mathrm{an}}$, equipped with induced topology
\[\{(x, e^{\vee}(x))\in Tot(L^{\vee})^{\mathrm{an}}: \abs{e^{\vee}(s)}(x) \leq \abs{s}_{\phi}(x)\cdot \mathrm{e}^{\epsilon}\ (resp. <\mathrm{e}^{\epsilon})\}\]
is called the \emph{dual closed (resp. open) disc bundle of radius $\mathrm{e}^{\epsilon}$} of the pseudometrized pair $(L,\phi)$, where $s$ is a local section of $L$. We denote it by $\overline{\mathbb{D}}^{\vee}(L, \phi, \epsilon)$ (resp. $\mathbb{D}^{\vee}(L, \phi, \epsilon))$.

\item Let $f:Y\rightarrow X$ be a morphism of separated $k$-schemes of finite type. Let $L$ be an invertible $\mathscr{O}_X$-module, equipped with a pseudometric $\phi$. We define a pseudometric $f^*\phi$ on $f^*(L)$ such that, for any section $s$ of $L$ on a Zariski open subset $U$ of $X$, one has
\[\forall\,y\in f^{-1}(U)^{\mathrm{an}},\quad  |f^*(s)|_{f^*\phi}(y)=|s|_{\phi}( f^{\mathrm{an}}(y)).\]
Since $f^{\mathrm{an}}:Y^{\mathrm{an}}\rightarrow X^{\mathrm{an}}$ is continuous (Proposition \ref{morphism analytification}), if the metric $\phi$ is continuous, so is $f^*\phi$. If $Y$ is a subscheme of $X$ and if $f:Y\rightarrow X$ is the canonical immersion, the restricted metric $f^*\phi$ is also denoted by $\phi|_Y$.  

\item\label{Item: metric and weight} Any map $f:X^{\mathrm{an}}\rightarrow\mathbb R\cup\{+\infty\}$ determines a pseudometric $\tau_f$ on $\mathscr{O}_X$ such that, for any regular function $a$ of $X$ on a Zariski open subset $U$, one has (with the convention $\mathrm{e}^{-\infty}=0$)
\[\forall\,x\in U^{\mathrm{an}},\quad \abs{a}_{\phi_f}(x)=|a|(x)\cdot \mathrm{e}^{-f(x)}.\]  Note that $f \mapsto\tau_f$ defines a bijection between the set of maps $X^{\mathrm{an}}\rightarrow \mathbb R\cup\{+\infty\}$ and that of pseudometrics on $\mathscr{O}_X$, which maps the set of real-valued functions bijectively to that of pseudometrics on $\mathscr{O}_X$. Moreover, a pseudometric $\phi_f$ is continuous if and only if $f$ is  continuous on $X^{\mathrm{an}}$. The trivial invertible sheaf $\mathscr{O}_X$ equipped with the pseudometric $\tau_f$ is denoted by $\mathscr{O}_X(f)$. The metric corresponding to the identically vanishing function is called the \emph{trivial metric} on $\mathscr{O}_X$.

\item\label{Item: distance between metrics} Let $\phi_1$ and $\phi_2$ be two metrics on $L$. The \emph{distance} of these two pseudometrics is a generalized positive real number (in $\R_+\cup \parenth{+\infty}$) defined by
\[\dist (\phi_1, \phi_2)=\sup_{x\in X^{\mathrm{an}}} \Big|\log \Big|\frac{\phi_1(x)}{\phi_2(x)}\Big|_{\widehat{\kappa}(x)}\Big|.\]
If $X$ is proper and $\phi_1$, $\phi_2$ are continuous metrics, then $\dist (\phi_1, \phi_2)\in \R_+$.
 
\item\label{addition of metrics} Let $L_1$ and $L_2$ be invertible $\mathscr{O}_X$-modules, and $\phi_1$ and $\phi_2$ be pseudometrics on $L_1$ and $L_2$ respectively. The pseudometric $\phi_1$ and $\phi_2$ induce by passing to tensor product a metric on $L_1\otimes L_2$, denoted by $\phi_1+\phi_2$. For any Zariski open subset $U$ of $X$ and any $(s_1,s_2)\in\Gamma(U,L_1)\times\Gamma(U,L_2)$, one has 
\[\forall\,x\in U^{\mathrm{an}},\quad\abs{s_1\cdot s_2}_{\phi_1+\phi_2}(x)=\abs{s_1}_{\phi_1}(x)\cdot\abs{s_2}_{\phi_2}(x).\]

If $\phi_1$ and $\phi_2$ are continuous, then $\phi_1+\phi_2$ is also continuous. 

In particular, for any $\epsilon\in \R$, we denote by $\phi(\epsilon)$ the pseudometric $\phi+\tau_{\mathrm{e}^{\epsilon}}$ on $L\otimes \mathscr{O}_X=L$.

Moreover, any metric $\phi$ on $L$ determines by passing to its dual a metric $-\phi$ on $L^{\vee}$ such that, for any Zariski open subset $U$ of $X$ and any $(s,\alpha)\in\Gamma(U,L)\times\Gamma(U,L^\vee)$, one has
\[\forall\,x\in U^{\mathrm{an}},\quad \abs{\alpha(s)}(x)=\abs{\alpha}_{-\phi}(x)\cdot\abs{s}_\phi(x).\]
If the metric $\phi$ is continuous, so is $-\phi$.  

\item\label{power and root of metric} 
Let $L$ be an invertible $\mathscr{O}_X$-module, and $n\in\mathbb N\setminus\{0\}$. A pseudometric $\phi$ on $L$ determines by tensor power a pseudometric on $L^{\otimes n}$ for any $n\in\mathbb N\setminus\{0\}$, denoted by $n\phi$. By convention, $0\phi$ denotes the trivial metric on $L^{\otimes 0}\cong\mathscr{O}_X$ (see \ref{Item: metric and weight} above).

Similarly, assume given a pseudometric $\phi$ on $L^{\otimes n}$. We denote by $\frac 1n\phi$ the pseudometric on $L$ such that, for any Zariski open subset $U$ of $X$ and any section $s\in\Gamma(U,L)$, one has
\[\norm{s}_{\frac 1n\phi}=\norm{s^n}_{\phi}^{1/n}.\]
If the pseudometric $\phi$ is continuous, then also is $\frac 1n\phi$.

\item\label{localization to fiber}
Let $L$ be an invertible $\mathscr{O}_X$-module. For any  $n$ such that $L^{\otimes n}$ is globally generated, let $\norm{\ndot}_{n}$ be a norm on $V_n(L)$. For any $x\in X^{\mathrm{an}}$, the evaluation map
\[V_n(L)\otimes_k \widehat{\kappa}(x)\longrightarrow L^{\otimes n}(x)\]
induces a quotient norm of $\norm{\ndot}_{n, \widehat{\kappa}(x)}$ on the $\widehat{\kappa}(x)$-vector space $L^{\otimes n}(x)$, denoted by $\norm{\ndot}_{n, X|x}$. This gives rise to a metric on $L^{\otimes n}$, which we call the \emph{Fubini-Study metric associated with $\norm{\ndot}_{n}$} on $L^{\otimes n}$, denoted by $\mathrm{FS}(\norm{\ndot}_{n})(x)$. The metric $\frac{1}{n}\mathrm{FS}(\norm{\ndot}_{n})$ on $L$ is called the \emph{$n$-th Fubini-Study metric associated with $\norm{\ndot}_{n}$} on $L$.

In particular, let $\parenth{s_{n, j}}_{j\in \Card{d_n}}$ be a basis of $V_n(L)$ and let $\norm{\ndot}_n$ be a ultrametric norm on $V_1(L)$ with respect to which $\parenth{s_{1, j}}_{j\in \Card{d_1}}$ is an orthogonal basis. Such a metric $\mathrm{FS}(\norm{\ndot}_{n})$ is studied and heavily used in \cite{CM}, and is said to be \emph{diagonalizable} in \cite{BE}. 

\item\label{FS envelop metric}
Similarly, let $\algnorm{\ndot}$ be an algebra norm on $V_{\sbullet}(L)$. For any $x\in X^{\mathrm{an}}$, the evaluation map induces a $\widehat{\kappa}(x)$-algebra homomorphism
\[V_{\sbullet}(L)\otimes \widehat{\kappa}(x)\rightarrow \bigoplus_{n\in \N}L^{\otimes n}(x)=:V_{\sbullet}(L)(x)\]
This algebra homomorphism induces a quotient algebra norm of the scalar extension $\algnorm{\ndot}_{\widehat{\kappa}(x)}$ on $V_{\sbullet}(L)(x)$, denoted by  $\algnorm{\ndot}_{X|x}$. Let $\widehat{V}(L, \algnorm{\ndot})(x)$ denote the separated completion of $(V_{\sbullet}(L)(x), \algnorm{\ndot}_{X|x})$. Once a non-zero element $e_1(x)\in L(x)$ is chosen, the second algebra can be identified with $\widehat{\kappa}(x)[T]$ by sending $e_1(x)$ to $T$. 

\item\label{Item: FS metric and FS envelop metric}
Let $L$ be an invertible $\mathscr{O}_X$ module. Let $\{\norm{\ndot}_{n}\}_{n\in \N}$ be a familly of norms on $\{V_n(L)\}_{n\in \N}$. If the sequence of metrics $\{\frac{1}{n}\mathrm{FS}(\norm{\ndot}_{n})\}_{n\in \N}$ converges pointwisely to a limit metric, we denote it by $\mathcal{P}(\{\norm{\ndot}_{n}\}_{n\in \N})$ and call it the \emph{Fubini-Study envelop metric} associated with $\{\norm{\ndot}_{n}\}_{n\in \N}$. 

Note that if the convergence is uniform for $x\in X^{\mathrm{an}}$, since Fubini-Study metrics $\{\frac{1}{n}\mathrm{FS}(\norm{\ndot}_{n})\}_{n\in \N}$ are continuous, the envelop metric will also be continuous. Conversely, if $X$ is proper over $\spec k$ and the envelop metric is continuous, then the convergence is uniform in $x\in X^{\mathrm{an}}$ as $X^{\mathrm{an}}$ is Hausdorff and compact by Theorem \ref{Analytification Hausdorff compact}. A metric $\phi$ on $L$ is \emph{asymptotic Fubini-Study} if it is a Fubini-Study envelop metric and the convergence is uniform for $x\in X^{\mathrm{an}}$ (see \cite[Definition 6.1]{BE}). Asymptotic Fubini-Study metrics are thus continuous. Note that asymptotic Fubini-Study property in this sense is equivalent to the notion of \emph{semipositive} metric by the terminology of \cite{CM}. We refer to \cite[\S 5.4]{BFJ} and \cite[\S 6.1]{BE} for a clear discussion of other various notions of semipositivity that have been proposed and studied in \cite{Zhang}, \cite{Gub}, \cite{Mor}, \cite{BFJ}, \cite{CLD}, \cite{BMPS}, \cite{CM}, \cite{GM} and literature therein.

In particular, let $\algnorm{\ndot}$ be an algebra seminorm on $V_{\sbullet}(L)$, and let $\{\norm{\ndot}_{n}\}_{n\in \N}$ be the associated familly of seminorms on $\{V_n(L)\}_{n\in \N}$. The seminorms $\{\frac{1}{n}\mathrm{FS}(\norm{\ndot}_{n})(x)\}_{n\in \N}$ on $L(x)$ satisfy sub-multiplicative property, so they converges to a limit seminorm on $L(x)$. This gives rise to a pseuodometric on $L$, called the \emph{Fubini-Study envelop pseudometric associated with $\algnorm{\ndot}$}. We denote it by $\mathcal{P}(\algnorm{\ndot})$. It is not necessarily continuous.

\item\label{Item: sup norm} Assume that $X$ is proper over $\spec k$. Note that $X^{\mathrm{an}}$ is then a compact Hausdorff space (see \cite[Theorem 3.4.8]{Ber}). Let $L$ be an invertible $\mathscr{O}_X$-module and $\phi$ be an upper semicontinuous metric on $L$ (see \ref{Item: metric} above). As $X^{\mathrm{an}}$ is compact, any upper semicontinuous function on $X^{\mathrm{an}}$ is bounded from above and attains its maximal value. In particular, for any $s\in V_1(L)$, one has
\[\norm{s}_{\phi}:=\sup_{x\in X^{\mathrm{an}}}\abs{s}_\phi(x)<+\infty.\] 
Moreover, $\norm{\ndot}_{\phi}:V_1(L)\rightarrow\mathbb R_{\geq 0}$ is a norm on $V_1(L)$. This norm is ultrametric since the absolute value $\abs{\ndot}$ on $k$ is non-Archimedean and $L$ is of rank $1$. We denote by $\algnorm{\ndot}_{\phi}$ the norm on the $k$-vector space $V_\sbullet(L)$ defined as 
\[\forall\,\underline{s}=(s_n)_{n\in\mathbb N}\in V_\sbullet(L),\quad\algnorm{\underline{s}}_{\phi}:=\sup_{n\in\mathbb N}\norm{s_n}_{n\phi}.\]
Note that the $k$-algebra $V_\sbullet(L)$ equipped with this norm forms a normed $k$-algebra. In fact, since $0\phi$ is the trivial metric on $\mathscr{O}_X$, one has $\algnorm{\mathbf{1}}_{0\phi}=1$, where $\mathbf{1}$ denotes the unit section of $\mathscr{O}_X$. Moreover, for $s_n\in V_n(L)$ and $s_m\in V_m(L)$, we have
\begin{equation*}
    \begin{split}
        \norm{s_m\cdot s_n}_{n\phi}=&\sup_{x\in X^{\mathrm{an}}}\abs{s_m}_{m\phi}(x)\cdot \abs{s_m}_{m\phi}(x)\\
        \leqslant& \sup_{x\in X^{\mathrm{an}}}\abs{s_m}_{m\phi}(x)\cdot \sup_{x\in X^{\mathrm{an}}}\abs{s_n}_{n\phi}(x)=\norm{s_m}_{m\phi}\cdot \norm{s_n}_{n\phi}
    \end{split}
\end{equation*}
Then $\algnorm{\ndot}_{\phi}$ is an algebra norm by \ref{Item: algebra norm}

In addition, the familly of norms $\{\norm{\ndot}_{n\phi}\}_{n\in \N}$ satisfies the power-multiplicative property for homogeneous elements: 
\[\forall s_n\in V_n(L), \forall m\in \N, \quad \norm{(s_n)^m}_{nm\phi}=(\norm{s_n}_{n\phi})^m.\]
In fact, the algebra norm $\algnorm{\ndot}_{\phi}$ is power-multiplicative also for non-homogeneous elements (see Proposition \ref{sup algebra norm is PM}).

Moreover, for every $n\in \N$, there is a metric $\frac{1}{n}\mathrm{FS}(\norm{\ndot}_{n\phi})$ on $L$, namely the $n$-th Fubini-Study metrics  on $L$ associated with $\norm{\ndot}_{n\phi}$.

One denotes by $\widehat{V}_\sbullet(L,\phi)$ the separated completion of the normed $k$-algebra $(V_\sbullet(L),\algnorm{\ndot}_{\phi})$. More generally, if $V_\sbullet$ is a graded sub-$k$-algebra of $V_\sbullet(L)$, by abuse of notation we still denote by $\algnorm{\ndot}_{\phi}$ the restriction of the norm $\algnorm{\ndot}_{\phi}$ on $V_\sbullet$ and denote by $\widehat{V}_{\sbullet}(\phi)$ the separated completion of the normed algebra $(V_\sbullet,\algnorm{\ndot}_{\phi})$. The restricted norm is also power-multiplicative. 

In particular, for any $N\in \N\setminus\{0\}$, if one takes $V_{\sbullet}$ to be $\bigoplus_{n\in \N}V_{nN}(L)$, denoted by $V_{\sbullet}^{(N)}(L)$, we denote by $\widehat{V}_{\sbullet}^{(N)}(L, \phi)$ the separated completion of $(V_{\sbullet}^{(N)}(L), \algnorm{\ndot}_{\phi})$.

\item\label{Item: quotient norm}
Assume that $X$ is proper over $\spec k$. Let $L$ be an invertible $\mathscr{O}_X$-module and $\phi$ be an upper semicontinuous metric on $L$. Let $f:Y\rightarrow X$ be a morphism of $k$-schemes of finite type. Let $\norm{\ndot}_{n\phi, X|Y}$ be the quotient norm of $\norm{\ndot}_{n\phi}$ on $V_n(L_{X|Y})$. It is ultrametric. Let $\algnorm{\ndot}_{\phi, X|Y}$ be the quotient algebra norm of $\algnorm{\ndot}_{\phi}$ on $V_{\sbullet}(L_{X|Y})$. In fact,
\[\forall \underline{t}=(t_n)_{n\in \N}\in V_{\sbullet}(L_{X|Y}), \quad  \algnorm{\underline{t}}_{\phi,X|Y}=\sup_{\begin{subarray}{c} n\in \N \end{subarray}}\norm{t_n}_{n\phi,X|Y}.\]
One denotes by $\widehat{V}_{\sbullet}(L_{X|Y}, \phi_{X|Y})$ the separated completion of the normed $k$-algebra $(V_{\sbullet}(L_{X|Y}), \algnorm{\ndot}_{\phi_{X|Y}})$. In particular, for any $N\in \N$, we denote by $V_{\sbullet}^{(N)}(L_{X|Y})$ the graded $k$-algebra $\bigoplus_{n\in \N}V_{nN}(L_{X|Y})$. This is a sub-algebra of $V_{\sbullet}(L_{X|Y})$. The restriction of $\algnorm{\ndot}_{\phi, X|Y}$ on this sub-algebra is still denoted by $\algnorm{\ndot}_{\phi, X|Y}$. One denotes by $\widehat{V}_{\sbullet}^{(N)}(L_{X|Y}, \phi_{X|Y})$ the separated completion of $(V_{\sbullet}^{(N)}(L_{X|Y}), \algnorm{\ndot}_{\phi, X|Y})$. In particular, if $f$ is the canonical immersion associated with a sub-scheme, we get a Banach $k$-algebra $\widehat{V}_{\sbullet}^{(N)}(L_{X|Y}, \phi_{X|Y})$.

\end{enumerate}

In the rest of the article, we make the following assumptions. For algebro-geometric data: let $X$ be an integral projective scheme over $\spec k$ of pure dimension $d$, $Y$ be a reduced closed sub-scheme of $X$ with its canonical closed immersion $i_Y: Y\rightarrow X$, and $L$ be an ample invertible $\mathscr{O}_X$-module. One can find $M\in \N$ such that $L^{\otimes M}$ is very ample and for any $n\geq M$, the restriction map from $V_n(L)$ to $V_n(L|_Y)$ is surjective, so $V_n(L|_Y)=V_n(L_{X|Y})$. For the metric data, let $\phi$ be an upper-semicontinuous metric on $L$.

\subsection{Algebraic properties of normed section algebra}\hfill\break
We show the power-multiplicativity of the supremum algebra norm $\algnorm{\ndot}_{\phi}$, and that the normed section algebras are reduced Banach algebras.

\begin{proposition}\label{sup algebra norm is PM}
    The algebra norm $\algnorm{\ndot}_{\phi}$ is power-multiplicative. Hence its spectral algebra seminorm is equal to itself, and it is an algbra norm.
\end{proposition}
\begin{proof}
    In fact, let $\underline{s}=(s_n)_{n\in \N}$ be an element of $V_{\sbullet}(L)$ and $m\in \N$, let $n_0\in \N$ be the smallest integer for which $\algnorm{\underline{s}}_{\phi}=\norm{s_{n_0}}_{n_0\phi}$. By the ultrametricity of $\algnorm{\ndot}_{\phi}$ and the power-multiplicativity of $\{\norm{\ndot}_{n\phi}\}_{n\in \N}$, one has
    \[\algnorm{\underline{s}^m}_{\phi}\leqslant\norm{(s_{n_0})^m}_{mn_0\phi}=(\norm{s_{n_0}}_{n_0\phi})^m=\algnorm{\underline{s}}_{\phi}^m.\]
    By the choice of $n_0$, one has
    \[\forall (j_0, \dots, j_l)\in \N^{l+1}, \sum_{i\in \Card{l}}i\cdot j_i=mn_0,\quad \Big\lVert\prod_{i\in \Card{l}}(s_i)^{j_i}\Big\rVert_{mn_0\phi}\leq \norm{(s_{n0})^m}_{mn_0\phi}\]
    and the equality holds if and only if $(j_0, \dots, j_l)=(0, \dots, 0, n_0, 0, \dots, 0)$ where $n_0$ is on the $m$-th place, so by the definition of $\algnorm{\ndot}$ and its ultra-metricity, one get
    \[\algnorm{\underline{s}^m}_{\phi}\geqslant \norm{(\underline{s}^m)_{mn_0}}_{mn_0\phi}\geq \norm{(s_{n0})^m}_{mn_0\phi}=\algnorm{\underline{s}}_{\phi}^m,\]
    hence there is an equality.
\end{proof}

\begin{corollary}\label{normed section algebra semi-simple}
    Then the Banach $k$-algebras  $\widehat{V}_{\sbullet}(L, \phi)$, $\widehat{V}_{\sbullet}(L|_Y, \phi|_Y)$ and $\widehat{V}_{\sbullet}(L_{X|Y}, \phi_{X|Y})$ are semi-simple. In particular, they are reduced.
\end{corollary}
\begin{proof}
    Let $\underline{s}\in \widehat{V}_{\sbullet}(L, \phi)$ be an element in $\rad(\widehat{V}_{\sbullet}(L, \phi))$, then by Proposition \ref{sup algebra norm is PM}, one has
    \[0=\algnorm{\underline{s}}_{\phi; \mathrm{sp}}=\algnorm{s}_{\phi},\]
    so 
    \[\forall n\in \N, \quad \norm{s_n}_{n\phi}=0.\]
    By the assumption, all componets $s_n$ are zero sections. So $\underline{s}=\underline{0}$, hence $\widehat{V}_{\sbullet}(L, \phi)$ is semi-simple. Same arguments works for $\widehat{V}_{\sbullet}(L|_Y, \phi|_Y)$.
    
    Let $\underline{t}\in \rad(\widehat{V}_{\sbullet}(L_{X|Y}, \phi_{X|Y}))$. Then $t_n\in  \rad(\widehat{V}_{\sbullet}(L_{X|Y}, \phi_{X|Y}))$ for every $n\in \N$. For any $m\in \N$, there exists $s_{nm}\in V_{nm}(L)$ such that $s_{nm}|_Y=t_n^m$ and 
    \[\lim_{\begin{subarray}{c}m\to \infty\end{subarray}}\norm{s_{nm}}_{nm\phi}^{\frac{1}{m}}=0.\]
    As $\norm{t_n^m}_{nm\phi|_Y}\leq \norm{s_{nm}}_{nm\phi}$, one has that \[\forall n\in \N, \quad\norm{t_n}_{n\phi|_Y}=0\]
    so $\underline{t}=\underline{0}$. Hence $\widehat{V}_{\sbullet}(L_{X|Y}, \phi_{X|Y})$ is semi-simple.
\end{proof}
\begin{remark}
    The reducity of closed sub-scheme $Y$ is necessary for the semi-simplicity of the Banach $k$-algebra $\widehat{V}_{\sbullet}(L_{X|Y}, \phi_{X|Y})$. 
\end{remark}

\subsection{Spectrum of normed section algebra}\hfill\break
We embed the Berkovich spectrum of normed section algebra into the analytification of the spectrum of the section algebra.

\begin{lemma}
    The homomorphism of inclusion of $k$-algebras $V(L)\rightarrow \widehat{V}_{\sbullet}(L,\phi)$ induces a continuous map between topological spaces $\M(\widehat{V}_{\sbullet}(L,\phi))\rightarrow (\spec V(L))^{\mathrm{an}}$ which is closed. Moreover, the map is injective.
\end{lemma}
\begin{proof}
    This is clear by Proposition \ref{normed spectrum into non-normed spectrum}.
\end{proof}

\begin{proposition}\label{general box all alg}
    There exist algebra norms $\algnorm{\ndot}_{\phi}^{\mathrm{aff}}$ on $V_{\sbullet}(L)$ and $\algnorm{\ndot}_{\phi, X|Y}^{\mathrm{aff}}$ on $V_{\sbullet}(L_{X|Y})$ such that the separated completions of
    $(V_{\sbullet}(L), \algnorm{\ndot}_{\phi}^{\mathrm{aff}})$ and $(V_{\sbullet}(L_{X|Y}), \algnorm{\ndot}_{\phi, X|Y}^{\mathrm{aff}})$ are affinoid algebras. Denote them by $\widehat{V}_{\sbullet}(L, \phi^{\mathrm{aff}})$ and $\widehat{V}_{\sbullet}(L_{X|Y}, \phi_{X|Y}^{\mathrm{aff}})$. Moreover, there exists a commutative diagram of homomorphisms of $k$-algebras with $\sigma$ and $\sigma_Y$ being homomorphisms of Banach $k$-algebras
    \[
    \xymatrix{
    V_{\sbullet}(L) \ar[d]^{i_Y} \ar[r] &\widehat{V}_{\sbullet}(L, \phi^{\mathrm{aff}}) \ar[r]^{\sigma} \ar[d]^{i_Y(\phi^\mathrm{aff})} &\widehat{V}_{\sbullet}(L, \phi) \ar[d]^{i_Y(\phi)} & \\
    V_{\sbullet}(L_{X|Y}) \ar[r] &\widehat{V}_{\sbullet}(L_{X|Y}, \phi_{X|Y}^{\mathrm{aff}}) \ar[r]^{\sigma|_Y} &\widehat{V}_{\sbullet}(L_{X|Y}, \phi_{X|Y})
    } 
    \]
    All homomorphisms (except $i_Y$) have dense images, and $i_Y$ is surjective.
\end{proposition}
\begin{proof}
    As $V_{\sbullet}(L)$ is a finitely generated sub-$k$-algebra of $\widehat{V}_{\sbullet}(L, \phi)$ which is dense for the topology induced by Banach algebra norm, by Proposition \ref{affinoid box}, there exists an affinoid algebra norm $\algnorm{\ndot}_{\phi}^{\mathrm{aff}}$ on $V_{\sbullet}(L)$ and a homomorphism of Banach algebras
    \[\sigma: \widehat{V}_{\sbullet}(L, \phi^{\mathrm{aff}}) \rightarrow \widehat{V}_{\sbullet}(L, \phi).\]
    
    One then takes $\algnorm{\ndot}_{\phi, X|Y}^{\mathrm{aff}}$ on $V_{\sbullet}(L_{X|Y})$ to be the quotient norm of $\algnorm{\ndot}_{\phi}^{\mathrm{aff}}$. By Example \ref{quotient of affinoid is affinoid}, this quotient norm is also an affinoid algebra norm. By this quotient construction, there exists a homomorphism of Banach $k$-algebras
    \[\sigma|_Y: \widehat{V}_{\sbullet}(L_{X|Y}, \phi_{X|Y}^{\mathrm{aff}}) \rightarrow\widehat{V}_{\sbullet}(L_{X|Y}, \phi_{X|Y})\]
    which fits into a commutative diagram with $i_Y(\phi)$ and $i_Y(\phi^\mathrm{aff})$.
\end{proof}

\begin{corollary}\label{general box all geom}
    The commutative diagram of homomorphisms of $k$-algebras induces a commutative diagram of continuous maps of topological spaces
    \[
    \xymatrix{
    (\spec(V_{\sbullet}(L)))^{\mathrm{an}} 
    &\M(\widehat{V}(L, \phi^{\mathrm{aff}})) \ar[l] &\M(\widehat{V}_{\sbullet}(L, \phi)) \ar[l]_{\sigma^*}\\
    (\spec(V_{\sbullet}(L_{X|Y})))^{\mathrm{an}} \ar[u]^{i_Y^*} &\M(\widehat{V}_{\sbullet}(L_{X|Y}, \phi_{X|Y}^{\mathrm{aff}})) \ar[l] \ar[u]^{i_Y(\phi^\mathrm{aff})^*} &\M(\widehat{V}_{\sbullet}(L_{X|Y}, \phi_{X|Y})) \ar[l]_{\sigma|_Y^*} \ar[u]^{i_Y(\phi)^*}
    } 
    \]
    All maps are closed. If the algebra seminorm $\algnorm{\ndot}_{\phi}$ is a norm, then all maps are injective.
\end{corollary}
\begin{proof}
    This follows from Proposition \ref{morphism of spectrum}, Proposition \ref{normed spectrum into non-normed spectrum} and Proposition \ref{injective}.
\end{proof}

\subsection{Fubini-Study metrics}\hfill\break
We study distances between Fubini-Study metrics, and gives explicit expressions for Fubini-Study metrics admitting non-Archimedean orthogonal basis. Many of the results here are also obtained in \cite{CM} or in \cite[Section 6]{BE}. 
\begin{lemma}\label{rk 1 field change}
    Assume that $L$ is globally generated. Let $\norm{\ndot}_1$ be a norm on $V_1(L)$ and let $\mathrm{FS}(\norm{\ndot}_{1})$ be the associated Fubini-Study metric. Then for any $x\in X^{\mathrm{an}}$ and $e(x)\in L(x)\setminus 0$, 
    \[\abs{e(x)}_{\mathrm{FS}(\norm{\ndot}_{1})}=\inf_{\begin{subarray}{c}\lambda\in \widehat{\kappa}(x),\ s_1\in V_1(L)\\ s_1(x)=\lambda\cdot e(x) \end{subarray}}\abs{\lambda}^{-1}\cdot\norm{s_1}.\]
    (with the convention that $0^{-1}=+\infty$)
\end{lemma}
\begin{proof}
    This follows from Lemma \ref{rk 1 quotient norm}.
\end{proof}

\begin{lemma}\label{FS is continuous}
    Let $\norm{\ndot}_n$ be a norm on $V_n(L)$, then $\frac{1}{n}\mathrm{FS}(\norm{\ndot}_{n})$ is a continuous metric on $L$. (\cite[Proposition 3.1]{CM})
\end{lemma}

\begin{lemma}\label{sup norm contractive}
    Let $\phi_1$ and $\phi_2$ be two upper-semicontinuous metrics on $L$, then 
    \[\dist (\norm{\ndot}_{\phi_1}, \norm{\ndot}_{\phi_2})\leq \dist (\phi_1, \phi_2).\]
    (see also \cite[Lemma 6.10]{BE})
\end{lemma}
\begin{proof}
    For any $s\in V_1(L)$, if $\norm{s}_{\phi_1}\geq \norm{s}_{\phi_2}$, let $x_1\in X^{\mathrm{an}}$ be a point where $\norm{s}_{\phi_1}$ is attained. One has
    \[\Big|\log \frac{\norm{s}_{\phi_1}}{\norm{s}_{\phi_2}}\Big| \leq \Big|\log\Big|\frac{s(x_1)}{s(x_1)}\Big|_{\widehat{\kappa}(x_1)}\Big| \leq\dist (\phi_1, \phi_2).\]
    Otherwise, let $x_2\in X^{\mathrm{an}}$ be a point where $\norm{s}_{\phi_2}$ is attained. Then
    \[\Big|\log \frac{\norm{s}_{\phi_2}}{\norm{s}_{\phi_1}}\Big|\leq \Big|\log\Big|\frac{s(x_2)}{s(x_2)}\Big|_{\widehat{\kappa}(x_2)}\Big| \leq\dist (\phi_1, \phi_2).\]
    Hence the desired inequality holds.
\end{proof}

\begin{lemma}\label{FS contractive}
    Let $\norm{\ndot}$ and $\norm{\ndot}'$ be two norms on $V_1(L)$, then
    \[\dist (\mathrm{FS}(\norm{\ndot}), \mathrm{FS}(\norm{\ndot}'))\leq \dist (\norm{\ndot}, \norm{\ndot}').\]
    (see also \cite[Equation (6.2)]{BE})
\end{lemma}
\begin{proof}
    For any $x\in X^{\mathrm{an}}$, let $e(x)\in L^{\mathrm{an}}(x)\setminus \parenth{0}$. Let $s, s'\in V_1(L)$ and $\lambda, \lambda'\in \hat{\kappa}(x)$ be elements such that
    \[s(x)=\lambda\cdot e(x),\  \norm{e(x)}_{X|x}=\abs{\lambda}^{-1}
    \norm{s},\] 
    \[s'(x)=\lambda'\cdot e(x),\  \norm{e(x)}_{X|x}=\abs{\lambda'}^{-1}\norm{s'}'.\]
    If $\norm{e(x)}_{X|x}>\norm{e(x)}'_{X|x}$, one has
    \begin{equation*}
        \begin{split}
            \dist(\norm{e(x)}_{X|x}, \norm{e(x)}'_{X|x})&=\Big|\log\frac{ \abs{\lambda}^{-1}\norm{s}}{\abs{\lambda'}^{-1} \norm{s'}'}\Big|
            \\&\leq \Big|\log\frac{\abs{\lambda}^{-1} \norm{s}}{\abs{\lambda}^{-1} \norm{s}'}\Big|\leq \dist(\norm{\ndot}, \norm{\ndot}').
        \end{split}
    \end{equation*}
    Otherwise, one has
    \begin{equation*}
        \begin{split}
            \dist(\norm{e(x)}_{X|x}, \norm{e(x)}'_{X|x})&=\Big|\log\frac{ \abs{\lambda}^{-1}\norm{s}}{\abs{\lambda'}^{-1} \norm{s'}'}\Big|
            \\&\leq \Big|\log\frac{\abs{\lambda'}^{-1} \norm{s'}}{\abs{\lambda'}^{-1} \norm{s'}'}\Big|\leq \dist(\norm{\ndot}, \norm{\ndot}').
        \end{split}
    \end{equation*}
    Varying $x$ and taking the supremum, one gets the desired inequality.
\end{proof}

\begin{proposition}\label{FS idempotent}
    Assume that there exist a norm $\norm{\ndot}_1$ on $V_1(L)$ such that $\phi$ is equal to $\mathrm{FS}(\norm{\ndot}_1)$. Then for any $n\in\N$,  $\mathrm{FS}(\norm{\ndot}_{n\phi})$ is equal to $n\phi$.
    (\cite[Proposition 3.3]{CM})
\end{proposition}

\begin{proposition}\label{FS envelop idempotent}
    Assume that $\phi$ is an asymptotic Fubini-Study metric on $L$, then the envelop metric $\mathcal{P}(\algnorm{\ndot}_{\phi})$ is equal to $\phi$.
    (see also \cite[Theorem 6.15 (iii)]{BE})
\end{proposition}
\begin{proof}
    By assumption, there exists a familly of norms $\{\norm{\ndot}_n\}_{n\in \N}$ such that uniformly for $x\in X^{\mathrm{an}}$,
    \[\quad\lim_{\begin{subarray}{c}n\to \infty\end{subarray}}\frac{1}{n}\mathrm{FS}(\norm{\ndot}_n)(\ast)(x)=\abs{\ast}_{\phi}(x),\]
    so for any $\epsilon>0$, there exists $N_0\in \N$ such that for any $n\geq N_0$
    \[\dist (n\phi, \mathrm{FS}(\norm{\ndot}_n))\leq n\epsilon,\]
    hence by Lemma \ref{sup norm contractive} and \ref{FS contractive},
    \[\dist (\mathrm{FS}(\norm{\ndot}_{n\phi}), \mathrm{FS}(\norm{\ndot}_{\mathrm{FS}(\norm{\ndot}_n)}))\leq n\epsilon.\]
    By Proposition \ref{FS idempotent}, one has $\mathrm{FS}(\norm{\ndot}_{\mathrm{FS}(\norm{\ndot}_n)}))=\mathrm{FS}(\norm{\ndot}_n)$, so
    \[\dist (\mathrm{FS}(\norm{\ndot}_{n\phi}), \mathrm{FS}(\norm{\ndot}_n))\leq n\epsilon, \]
    and
    \[\dist (\frac{1}{n}\mathrm{FS}(\norm{\ndot}_{n\phi}), \frac{1}{n}\mathrm{FS}(\norm{\ndot}_n))\leq \epsilon.\]
    Taking limit for $n\to \infty$ and then for $\epsilon\to 0$, one has 
    \[\dist (\mathcal{P}(\algnorm{\ndot}_{\phi}, \phi)=0,\]
    so the two metrics are equal.
\end{proof}

If the ultrametric norm $\norm{\ndot}_n$ admits an orthogonal basis, one can calculate explicitly the associated Fubini-Study metric $\mathrm{FS}(\norm{\ndot}_n)$.
\begin{lemma}\label{ultrametric ortho min max}
    Let $(K, \abs{\ndot}_K)$ be a complete ultrametric valued field extension of $(k, \abs{\ndot})$. Then for any $\parenth{r_j}_{j\in \Card{d}}$ elements in $\R_{+}$, one has
    \[\inf_{\sum_{j\in \Card{d}}\kappa_j=1}\max_{j}\Big\{\abs{\kappa_j}_{K}\cdot r_j\Big\}=\min_{j\in \Card{d}}\parenth{r_j},\]
    where $\parenth{\kappa_j}_{j\in \Card{d}}$ are elements in $K$.
\end{lemma}
\begin{proof}
    On the one hand, let $j_0\in \Card{d}$ be an index such that $r_{j_0}$ is minimal. By taking $\kappa_j=0$ for $j\neq j_{0}$ and $\kappa_{j_0}=1$, one sees that 
    \[\inf_{\sum_{j\in \Card{d}}\kappa_j=1}\max_{j}\Big\{\abs{\kappa_j}_{K}\cdot r_j\Big\}\leq r_{j_0}=\min_{j\in \Card{d}}\parenth{r_j}.\]
    On the other hand, by the ultrametricity of $\abs{\ndot}_K$, if $\sum_{j}\kappa_j=1$, then there exist at least one $j_1\in \Card{d}$ such that $\abs{\kappa_{j_1}}\geq 1$, so
    \begin{equation*}
        \begin{split}
            \inf_{\sum_{j\in \Card{d}}\kappa_j=1}\max_{j}\Big\{\abs{\kappa_j}_{K}\cdot r_j\Big\}&\geq \inf_{\sum_{j\in \Card{d}}\kappa_j=1}\abs{\kappa_{j_1}}_{K}\cdot r_{j_1}\\
            &\geq \inf_{\sum_{j\in \Card{d}}\kappa_j=1} r_{j_1}
            =r_{j_1}\geq\min_{j\in \Card{d}}\parenth{r_j}.
        \end{split}
    \end{equation*}
    Hence the two sides are equal.
\end{proof}
\begin{proposition}\label{orthogonal FS metric}
    Assume that $L^{\otimes n}$ is globally generated. Let $\parenth{s_{n,j}}_{j\in\Card{d_n}}$ be a basis of $V_n(L)$. Let $\norm{\ndot}_n$ be a ultrametric norm on $V_n(L)$ with respect to which $\parenth{s_{n,j}}_{j\in\Card{d_n}}$ is orthogonal. Then for any $x\in X^{\mathrm{an}}$ and $e_n(x)\in L^{\otimes n}(x)\setminus 0$,
    \[\abs{e_n(x)}_{\mathrm{FS}(\norm{\ndot}_n)}=\min_{j\in\Card{d_n}}\Big\{\Big|\frac{e_n(x)}{s_{n,j}(x)}\Big|_{\widehat{\kappa}(x)}\cdot \norm{s_{n, j}}_n \Big\},\]
    with the convention that $0^{-1}=+\infty$.
    (see also \cite[Lemma 3.3]{CM})
\end{proposition}
\begin{proof}
    For $j\in\Card{d_n}$, let $\lambda_j\in\widehat{\kappa}(x)$ be such that $s_{n,j}(x)=\lambda_j\cdot e_n(x)$. By construction,
    \begin{equation*}
        \begin{split}
            \abs{e_n(x)}_{\mathrm{FS}(\norm{\ndot}_n)}
            =&\inf_{\mu_j\in \widehat{\kappa}(x), \ \sum\mu_j\cdot s_{n,j}(x)=e_n(x)}\Big\|\sum_{j\in \Card{d_n}}\mu_j\cdot s_{n,j}\Big\| \\
            =&\inf_{\mu_j\in \widehat{\kappa}(x), \ \sum\mu_j\cdot s_{n,j}(x)=e_n(x)}\max_{j}\Big\{\abs{\mu_j}_{\widehat{\kappa}(x)}\cdot \norm{s_{n,j}}\Big\}\\
            =&\inf_{\kappa_j\in \widehat{\kappa}(x), \ \sum\kappa_j=1}\max_{j}\Big\{\abs{\kappa_j\cdot\lambda_j^{-1}}_{\widehat{\kappa}(x)}\cdot \norm{s_{n,j}}\Big\}\\
            =&\inf_{\kappa_j\in \widehat{\kappa}(x), \ \sum\kappa_j=1}\max_{j}\Big\{\abs{\kappa_j}_{\widehat{\kappa}(x)}\cdot \abs{\lambda_j}_{\widehat{\kappa}(x)}^{-1}\cdot \norm{s_{n,j}}\Big\}\\
            =&\min_{j\in\Card{d_n}}\Big\{\abs{\lambda_j}_{\widehat{\kappa}(x)}^{-1} \cdot \norm{s_{n, j}}_n \Big\}.
        \end{split}
    \end{equation*}
    The last equality is obtained with Lemma \ref{ultrametric ortho min max}.
\end{proof}
\begin{corollary}\label{orthogonal FS metric dual}
    With the same hypothesis as above, for $e_n^{\vee}(x)\in (L^{\otimes n})^{\vee}(x)\setminus 0$, one has
    \[\abs{e_n^{\vee}(x)}_{\mathrm{FS}(\norm{\ndot}_n)^{\vee}}=\max \Big\{\abs{e_n^{\vee}(x)(s_{n,j})}_{\widehat{\kappa}(x)}\cdot \norm{s_{n, j}}_n^{-1}\Big\}.\]
\end{corollary}
\begin{proof}
    It suffices to note that $e_n^{\vee}(x)(s_{n,j})=\lambda_j$.
\end{proof}

\subsection{Dual unit disc bundle}\hfill\break
Let $\algnorm{\ndot}$ be an algebra norm on $V_{\sbullet}(L)$, such that $V_n(L)$ are orthogonal subspaces for $\algnorm{\ndot}$. We relate the Berkovich spectrum of normed section algebra with the dual unit disc bundle with respect to the envelop metric.

\begin{proposition}
    Let $z\in (\spec V_{\sbullet}(L))^{\mathrm{an}}$ be a point. Let $(x, e^{\vee}(x))\in Tot(L^{\vee})$ be the point $(p(\pmb{0})^{-1})^{\mathrm{an}}(z)$. Then $z\in\M(\widehat{V}(L, \algnorm{\ndot}))$ if and only if one of the following criteria holds
    \begin{enumerate}
        \item \label{point in disc bundle global} there exist $C(z)>0$ such that
        \[\forall \underline{s}\in V_{\sbullet}(L), \quad \abs{\underline{s}}_z\leq C(z)\cdot \algnorm{\underline{s}}.\]
        \item \label{point in disc bundle local} there exist $C(z)>0$ such that
        \[\forall \underline{s}(x)\in V_{\sbullet}(L)(x),\quad\abs{\underline{s}(x)}_z\leq C(z)\cdot \algnorm{\underline{s}(x)}_{X|x}.\]
        \item \label{point in disc bundle metric} there exist $C'(z)=1$ such that 
        \[\forall e_1(x)\in V_{1}(L)(x),\quad\abs{e_1(x)}_z\leq \algnorm{e_1(x)}_{(X|x);\mathrm{sp}},\]
        where $\algnorm{\ndot}_{(X|x);\mathrm{sp}}$ is the spectral algebra seminorm of $\algnorm{\ndot}_{X|x}$.
    \end{enumerate}
\end{proposition}
\begin{proof}
    The criterion \ref{point in disc bundle global} unfolds the definition of the fact that $z\in \M$. The criterion \ref{point in disc bundle local} is equivalent to the criterion \ref{point in disc bundle global}, as $\algnorm{\ndot}_{X|x}$ is the quotient algebra norm of $\algnorm{\ndot}_{\widehat{\kappa}(x)}$ for the evaluation map $\mathrm{ev}(x)$. The criterion \ref{point in disc bundle metric} is equivalent to the criterion \ref{point in disc bundle local}: if \ref{point in disc bundle local} holds, then 
    \[\forall n\in\N, \quad \abs{e_1(x)}_z\leq C(z)^{\frac{1}{n}}\cdot\algnorm{e_1^{\otimes n}(x)}_{X|x}^{\frac{1}{n}},\]
    so \ref{point in disc bundle metric} holds after a limit process for $n\to \infty$. Conversely, if \ref{point in disc bundle metric} holds, then since $V_n(L)(x)$ is spaned by $e_1^{\otimes n}(x)$ over $\widehat{\kappa}(x)$, one has
    \[\forall n\in \N, \quad \abs{s_n}_z\leq \algnorm{s_n(x)}_{(X|x);\mathrm{sp}}\leq \algnorm{s_n(x)}_{X|x},\]
    so \ref{point in disc bundle local} holds by the ultra-metricity of $\abs{\ndot}_z$ and the orthogonality of $\algnorm{\ndot}_{X|x}$ for $V_n(L)$'s.
\end{proof}

\begin{corollary}
    With the same notations as above, the algebra seminorm $\algnorm{\ndot}_{(X|x);\mathrm{sp}}$ on $L(x)$ is equal to $\abs{\ndot}_{\mathcal{P}(\algnorm{\ndot})}(x)$.
\end{corollary}
\begin{proof}
    For any $x\in X^{\mathrm{an}}$ and any $e_1(x)\in V_1(L)(x)$, we have
    \begin{equation*}
        \begin{split}
            \algnorm{e_1(x)}_{(X|x);\mathrm{sp}}&=\lim_{\begin{subarray}{c}n\to \infty \end{subarray}}\algnorm{e_1^{\otimes n}(x)}_{X|x}^{\frac{1}{n}}\\
            &=\lim_{\begin{subarray}{c}n\to \infty \end{subarray}}\frac{1}{n}\mathrm{FS}(\norm{\ndot}_n)(e_1(x))(x)=\mathcal{P}(\algnorm{\ndot})(e_1(x))(x).
        \end{split}
    \end{equation*}
\end{proof}
\begin{remark}
    The resulting algebra seminorm $\abs{\ndot}_{\mathcal{P}(\algnorm{\ndot})}(x)$ gives rise to a pseudometric on $L^{\mathrm{an}}$.
\end{remark}

\begin{lemma}\label{blow up topology}
    The map $p(\mathbf{0})^{\mathrm{an}}$ induces a continuous map of topological spacecs
    \[p(\pmb{0})^{\mathrm{an}}: Tot(L^{\vee})^{\mathrm{an}}\rightarrow (\spec V_{\sbullet}(L))^{\mathrm{an}}.\]
    Moreover, it induces a homeomorphism
    \[p(\pmb{0})^{\mathrm{an}}: Tot(L^{\vee})^{\mathrm{an}}\setminus \mathbb{O}^{\mathrm{an}}\rightarrow (\spec V_{\sbullet}(L))^{\mathrm{an}}\setminus \pmb{0}^{\mathrm{an}}.\]
\end{lemma}
\begin{proof}
    By Proposition \ref{morphism analytification}, the morphism $p(\pmb{0})$ of schemes of finite type over $\spec k$ induces a continuous map betweeen the topological space of their analytification:
    \[p(\pmb{0})^{\mathrm{an}}: Tot(L^{\vee})^{\mathrm{an}}\rightarrow (\spec V_{\sbullet}(L))^{\mathrm{an}}.\]
    Since 
    \[p(\pmb{0}): Tot(L^{\vee})\setminus \mathbb{O}\rightarrow \spec (V_{\sbullet}(L))\setminus \pmb{0}\]
    is an isomorphism of schemes of finite type, its analytification induces a homeomorphism by Proposition \ref{morphism analytification}
    \[p(\pmb{0})^{\mathrm{an}}: Tot(L^{\vee})^{\mathrm{an}}\setminus \mathbb{O}^{\mathrm{an}}\rightarrow (\spec V_{\sbullet}(L))^{\mathrm{an}}\setminus \pmb{0}^{\mathrm{an}}.\]
\end{proof}

\begin{proposition}\label{spectrum is disc bundle}
    The map $p(\mathbf{0})^{\mathrm{an}}$ induces a continuous map of topological spacecs
    \[\overline{\mathbb{D}}^{\vee}(L, \mathcal{P}(\algnorm{\ndot}), 0)\rightarrow (\spec V_{\sbullet}(L))^{\mathrm{an}} \]
    which induces a homeomorphism between
    \[\overline{\mathbb{D}}^{\vee}(L, \mathcal{P}(\algnorm{\ndot}), 0)\setminus \mathbb{O}^{\mathrm{an}}\rightarrow (\M(\widehat{V}_{\sbullet}(L,\algnorm{\ndot})))\setminus \mathbf{0}^{\mathrm{an}}. \]
\end{proposition}
\begin{proof}
    Starting with the continuous map in Lemma \ref{blow up topology}, we can determine the pre-image of $\M(\widehat{V}_{\sbullet}(L,\algnorm{\ndot}))$: let $z\in (\spec V_{\sbullet}(L))^{\mathrm{an}}\setminus\pmb{0}^{\mathrm{an}}$ be a point and $(x, e^{\vee}(x))$ be its unique pre-image under $p(\pmb{0})^{\mathrm{an}}$, where $x\in X^{\mathrm{an}}$ and $e^{\vee}(x)\in L^{\vee}(x)$. By Lemma \ref{point in disc bundle metric}, if we fix a non-zero element $e_1(x)\in L(x)$, the point $z$ lies in $\M(\widehat{V}_{\sbullet}(L,\algnorm{\ndot}))$ if and only if
    \[\abs{e_1(x)}_z\leq \abs{e_1(x)}_{\mathcal{P}(\algnorm{\ndot})}(x).\]
    This condition is equivalent to
    \[\abs{e^{\vee}(e_1)(x)}\leq \abs{e_1(x)}_{\mathcal{P}(\algnorm{\ndot})}(x).\]
    Hence there exists a continuous surjective map
    \[p(\pmb{0})^{\mathrm{an}}: \overline{\mathbb{D}}^{\vee}(L, \mathcal{P}(\algnorm{\ndot}), 0)\rightarrow (\M(\widehat{V}_{\sbullet}(L,\algnorm{\ndot}))).\]
    If we remove $\mathbb{O}^{\mathrm{an}}$ and $\pmb{0}^{\mathrm{an}}$ from the domain and image, the restricted map is indeed a homeomorphism. 
\end{proof}


One can give a precise description of dual unit disc bundle for a Fubini-Study metric admitting orthogonal basis.
\begin{proposition}\label{unit disc bundle of FS basis}
    Let $n\in\N$ be an integer such that $L^{\otimes n}$ is globally generated. Let $\parenth{s_{n, j}}_{j\in \Card{d_n}}$ be a basis of $V_n(L)$. Let $\norm{\ndot}_n$ be a ultrametric norm on $V_n(L)$ with respect to which this basis is orthogonal. Let $(x, e_1^{\vee}(x))\in Tot(L^{\vee})$ and $z\in (\spec V_{\sbullet}(L))^{\mathrm{an}}$ be it image under $p(\pmb{0})^{\mathrm{an}}$, then $(x, e_1^{\vee}(x))\in \overline{\mathbb{D}}^{\vee}(L, \frac{1}{n}\mathrm{FS}(\norm{\ndot}_n))$ (resp.$\mathbb{D}^{\vee}(L, \frac{1}{n}\mathrm{FS}(\norm{\ndot}_n))$) if and only if
    \[\forall j\in \Card{d_n},\ \abs{s_{n, j}(z)}\leq \norm{s_{n, j}}_n \ (\text{resp.}<\norm{s_{n, j}}_n).\]
    In particular, the image of $\mathbb{D}^{\vee}(L, \frac{1}{n}\mathrm{FS}(\norm{\ndot}_n))$ under $p(\pmb{0})^{\mathrm{an}}$ is an open subset in $(\spec V_{\sbullet}(L))^{\mathrm{an}}$.
\end{proposition}
\begin{proof}
    The assertion is clear if $e_1(x)=0$. For $e_1(x)\neq0$, let $e_n(x)=e_1^{\otimes n}(x)$, note that
    \[\abs{e_1^{\vee}(x)}_{\frac{1}{n}\mathrm{FS}(\norm{\ndot}_n)^{\vee}}=(\abs{e_n^{\vee}(x)}_{\mathrm{FS}(\norm{\ndot}_n)^{\vee}})^{\frac{1}{n}}.\]
    By Corollary \ref{orthogonal FS metric dual}, one has
    \[\abs{e_n^{\vee}(x)}_{\mathrm{FS}(\norm{\ndot}_n)^{\vee}}=\max \Big\{\abs{e_n^{\vee}(x)(s_{n,j})}_{\widehat{\kappa}(x)}\cdot \norm{s_{n, j}}_n^{-1}\Big\}\]
    so
    \[\abs{e_1^{\vee}(x)}_{\frac{1}{n}\mathrm{FS}(\norm{\ndot}_n)^{\vee}}=\max \Big\{\abs{e_n^{\vee}(x)(s_{n,j})}_{\widehat{\kappa}(x)}\cdot \norm{s_{n, j}}_n^{-1}\Big\}^{\frac{1}{n}}.\]
    Tautologically, one has
    \[s_{n, j}(z)=(e_1^{\otimes n})^{\vee}(x)(s_{n,j})=e_n^{\vee}(x)(s_{n,j}),\]
    so the criterion holds. By these defining equations, it is easy to see that the image of the open dual unit disc bundle is an open set.
\end{proof}

\begin{corollary}\label{0 is in the open}
    Let $\phi$ be a continuous metric on $L$. Then the image of $\mathbb{D}^{\vee}(L, \phi)$ under $p(\pmb{0})^{\mathrm{an}}$ is an open subset of $(\spec V_{\sbullet}(L))^{\mathrm{an}}$.
\end{corollary}
\begin{proof}
    As $\phi$ is continuous, $\mathbb{D}^{\vee}(L, \phi)\setminus \mathbb{O}^{\mathrm{an}}$ is an open subset of $Tot(L^{\vee})^{\mathrm{an}}\setminus \mathbb{O}^{\mathrm{an}}$. By Lemma \ref{blow up topology}, under the map $p(\pmb{0})^{\mathrm{an}}$, the image of $\mathbb{D}^{\vee}(L, \phi)\setminus \mathbb{O}^{\mathrm{an}}$ is an open subset of $(\spec V_{\sbullet}(L))^{\mathrm{an}}\setminus \pmb{0}^{\mathrm{an}}$, so it is also an open subset of $(\spec V_{\sbullet}(L))^{\mathrm{an}}$. It suffices to treat $\pmb{0}^{\mathrm{an}}$ which is the image of $\mathbb{O}^{\mathrm{an}}$.
    
    As $L$ is ample, there exist $n\in \N$ such that $L^{\otimes n}$ is globally generated. Let $\parenth{s_{n,j}}_{j\in \Card{d_n}}$ be a basis and let $\psi$ be the Fubini-Study metric associated with some ultrametric norm $\norm{\ndot}_n$ for which this basis is orthogonal. As both $\frac{1}{n}\psi$ and $\phi$ are continuous and $X^{\mathrm{an}}$ is compact, there exist $\alpha\in \R$ such that 
    \[\forall x\in X^{\mathrm{an}},\ \frac{1}{n}\psi(\alpha)(x)\leq \phi(x),\]
    so 
    \[p(\pmb{0})^{\mathrm{an}}(\mathbb{D}^{\vee}(L, \frac{1}{n}\psi(\alpha)))\subseteq p(\pmb{0})^{\mathrm{an}}(\mathbb{D}^{\vee}(L, \phi)).\]
    the left hand side is an open subset of $(\spec V_{\sbullet}(L))^{\mathrm{an}}$ by Lemma \ref{unit disc bundle of FS basis}. Then
    \[p(\pmb{0})^{\mathrm{an}}(\mathbb{D}^{\vee}(L, \phi))=p(\pmb{0})^{\mathrm{an}}(\mathbb{D}^{\vee}(L, \phi)\setminus\mathbb{O}^{\mathrm{an}})\cup p(\pmb{0})^{\mathrm{an}}(\mathbb{D}^{\vee}(L, \frac{1}{n}\psi(\alpha)))\]
    is an open set in $(\spec V_{\sbullet}(L))^{\mathrm{an}}$.
\end{proof}

\begin{corollary}\label{disc inclusion}
    If $\mathcal{P}(\phi)$ is continous, then for any $\epsilon>0$, one has
    \[\M(\widehat{V}_{\sbullet}(L, \phi))\subseteq \Int^{\mathrm{top}}_{V_{\sbullet}}(\M(\widehat{V}_{\sbullet}(L, \phi(\epsilon)))),\]
    where $\Int^{\mathrm{top}}_{V_{\sbullet}}$ denotes the topological interior as subspace of $\spec(V_{\sbullet}(L))^{\mathrm{an}}$.
\end{corollary}
\begin{proof}
    By Proposition \ref{spectrum is disc bundle}, the left hand side $\M(\widehat{V}_{\sbullet}(L, \phi))$ is identified with $p(\pmb{0})^{\mathrm{an}}(\overline{\mathbb{D}}^{\vee}(L, \phi))$, which is contained in the open subset $p(\pmb{0})^{\mathrm{an}}(\mathbb{D}^{\vee}(L, \phi(\epsilon)))$. This open subset is contained in $p(\pmb{0})^{\mathrm{an}}(\overline{\mathbb{D}}^{\vee}(L, \phi(\epsilon)))$ which is identified with $\M(\widehat{V}_{\sbullet}(L, \phi(\epsilon)))$, hence this open subset is contained in the topological interior of the later, the right hand side.
\end{proof}

\subsection{Comparison of algebra norms}\hfill\break

We compare the quotient algebra norm and the supremum algebra norm on the restricted section algebra, and get directly a (non-uniform) extension theorem.
\begin{proposition}\label{envelop quot is envelop sup restricted}
    Let $\phi$ be an upper-semicontinuous metric on $L$. Then
    \[\mathcal{P}(\algnorm{\ndot}_{\phi})|_Y=\mathcal{P}(\algnorm{\ndot}_{\phi, X|Y}).\]
\end{proposition}
\begin{proof}
    By definition, for any $y\in Y^{\mathrm{an}}$, one has
    \[\mathcal{P}(\algnorm{\ndot}_{\phi})(y)=\lim_{\begin{subarray}{c}n\to\infty \end{subarray}}\frac{1}{n}\mathrm{FS}(\norm{\ndot}_{n\phi})(y),\]
    \[\mathcal{P}(\algnorm{\ndot}_{\phi, X|Y})(y)=\lim_{\begin{subarray}{c}n\to\infty \end{subarray}}\frac{1}{n}\mathrm{FS}(\norm{\ndot}_{n\phi, X|Y})(y).\]
    Since the $k$-linear map $V_n(L)\rightarrow V_n(L_{X|Y})$ is surjective for all large $n\in \N$, and $\norm{\ndot}_{n\phi, X|Y}$ is the quotient norm of $\norm{\ndot}_{n\phi}$, one has
    \[\mathrm{FS}(\norm{\ndot}_{n\phi})(y)=\mathrm{FS}(\norm{\ndot}_{n\phi, X|Y})(y).\]
    Hence the two envelop metrics are equal.
\end{proof}

\begin{lemma}\label{asymptotic FS restriction}
    Let $\phi$ be an asymptotic Fubini-Study metric on $L$. Then $\phi|_Y$ is an asymptotic Fubini-Study metric on $L|_Y$.
\end{lemma}
\begin{proof}
    Suppose that $\phi$ is the pointwise limit on $X^{\mathrm{an}}$ of $\parenth{\frac{1}{n}\mathrm{FS}(\norm{\ndot}_n)}$, where $\parenth{\norm{\ndot}_n}$ are norms on $V_n(L)$. Then $L|_Y$ is the pointwise limit of $\parenth{\frac{1}{n}\mathrm{FS}(\norm{\ndot}_{n, X|Y})}$ on $Y^{\mathrm{an}}$.
\end{proof}

\begin{proposition}\label{Envelop quot=Envelop sup}
    Let $\phi$ be a asymptotic Fubini-Study metric on $L$. Consider two algebra norms $\algnorm{\ndot}_{\phi|_Y}$ and $\algnorm{\ndot}_{\phi, X|Y}$ on $V_{\sbullet}(L_{X|Y})$. Then the three metrics are equal
    \[\mathcal{P}(\algnorm{\ndot}_{\phi, X|Y})=\mathcal{P}(\algnorm{\ndot}_{\phi|_Y})=\phi|_Y.\]
\end{proposition}
\begin{proof}
    By Proposition \ref{envelop quot is envelop sup restricted}, 
    \[\mathcal{P}(\algnorm{\ndot}_{\phi, X|Y})=\mathcal{P}(\algnorm{\ndot}_{\phi})|_Y=\phi|_Y.\]
    It suffices to show the second equality. By Lemma \ref{asymptotic FS restriction}, $\phi|_Y$ is an asymptotic Fubini-Study metric on $L|_Y$. By Proposition \ref{FS envelop idempotent}, 
    \[\mathcal{P}(\algnorm{\ndot}_{\phi|_Y})=\phi|_Y.\]
\end{proof}

\begin{corollary}\label{spectral of quot is sup}
    Let $\phi$ be a asymptotic Fubini-Study metric on $L$. Then on $V_{\sbullet}(L_{X|Y})$, the spectral algebra seminorm of $\algnorm{\ndot}_{\phi, X|Y}$ is equal to $\algnorm{\ndot}_{\phi|_Y}$. There exists a canonical homeomorphism
    \[\M(\widehat{V}_{\sbullet}(L_{X|Y}, \phi_{X|Y}))\simeq \M(\widehat{V}_{\sbullet}(L_{X|Y}, \phi|_Y)).\]
\end{corollary}
\begin{proof}
    By Theorem \ref{uniformization keeps spectrum}, one has a homeomorphism
    \[\M(\widehat{V}_{\sbullet}(L_{X|Y}, \algnorm{\ndot}_{\phi, X|Y;\mathrm{sp}}))\simeq\M(\widehat{V}_{\sbullet}(L_{X|Y}, \algnorm{\ndot}_{\phi,X|Y}))=\M(\widehat{V}_{\sbullet}(L_{X|Y}, \phi_{X|Y})),\]
    by Proposition \ref{Envelop quot=Envelop sup}, one has
    \[\M(\widehat{V}_{\sbullet}(L_{X|Y}, \algnorm{\ndot}_{\phi, X|Y;\mathrm{sp}}))\simeq\M(\widehat{V}_{\sbullet}(L_{X|Y}, \algnorm{\ndot}_{\phi|_Y}))=\M(\widehat{V}_{\sbullet}(L_{X|Y}, \phi|_Y)).\]
    By Proposition \ref{Gelfand uniform}, the two power-multiplicative algebra seminorms $\algnorm{\ndot}_{\phi|_Y}$ and $\algnorm{\ndot}_{\phi, X|Y;\mathrm{sp}}$ on $V_{\sbullet}(L_{X|Y})$ are equal since they are both supremum norms on the same spectrum.
\end{proof}

\begin{theorem}\label{weak QE}
    Let $\phi$ be an asymptotic Fubini-Study metric on $L$, then for any $\epsilon>0$, and any $t_{1}\in V_{1}(L|_Y)$, there exists $n_Y\in \N$ such that for any $n\geq n_Y$, there exists $s_{n}\in V_{n}(L)$ with $s_{n}|_Y=t_1^{\otimes n}$ and
    \[\norm{s_{n}}_{n\phi}\leq \mathrm{e}^{n\epsilon}\cdot (\norm{t_1}_{\phi|_Y})^n.\]
\end{theorem}
\begin{proof}
    For any $M\leq m\leq 2M-1$, we have $t_1^{\otimes m}\in V_m(L_{X|Y})$. By Corollary \ref{spectral of quot is sup}, for any $\epsilon>0$, there exists $N_m\in \N$ such that for $l\geq N_m$,
    \[\algnorm{t_1^{\otimes ml}}_{\phi, X|Y}^{\frac{1}{l}}/\algnorm{t_1^{\otimes m}}_{\phi|_Y}\leq \epsilon.\]
    It is easy to see that there exists $n_Y\in \N$ such that the set of integers
    \[\{ml: l\geq N_m, M\leq m\leq 2M-1\}\]
    contains a subset of form $\N-\Card{n_Y-1}$: the case $M=1$ is clear; if $M\geq 2$, the fact that $M$ and $M+1$ are coprime guarantees the existence of $n_Y$. Note that $\algnorm{\ndot}_{\phi|_Y}$ is power-multiplicative, so for any $n\geq n_Y$, there exists $s_n\in V_n(L)$ with $s_n|_Y=t_1^{\otimes n}$ such that
    \[\norm{s_{n}}_{n\phi}\leq \mathrm{e}^{n\epsilon}\cdot \norm{t_1^{\otimes n}}_{n\phi|_Y}=\mathrm{e}^{n\epsilon}\cdot (\norm{t_1}_{\phi|_Y})^n.\]
\end{proof}
\begin{remark}
    This result is first obtained in \cite{CM}, by using approximation of $\phi$ by model metrics. Here we give another proof. Note that a slight unsatisfactory point of this version of metric extension theorem is that the degree $n_Y$ depends \textit{a priori} on the choice of initial data, the restricted section $t_1$. We will remove this dependence in the following sections.
\end{remark}

\section{Geometric approximation}

Recall that we shall compare two algebra norms $\algnorm{\ndot}_{\phi, X|Y}$ and $\algnorm{\ndot}_{\phi|_Y}$ on the restricted section algebra $V_{\sbullet}(L_{X|Y})$, where
\[V_{\sbullet}(L_{X|Y})=\bigoplus_{n\in \N}\mathrm{Im}(V_n(L)\xrightarrow{|_Y}V_n(L|_Y)),\]
the second being the spectral algebra norm of the first by Proposition \ref{spectral of quot is sup}. Some intermediate Banach algebra (with a uniform algebra norm) $\mathcal{W}$ is needed in the comparison.

For any $\epsilon>0$, we begin approximating $\M(\widehat{V}_{\sbullet}(L_{X|Y}, \phi_{X|Y}))$ by a special domain $W_{\epsilon}$, in order to construct a homomorphism of Banach algebras from $\widehat{V}_{\sbullet}(L_{X|Y}, \phi_{X|Y})$ to $\mathcal{W}_{\epsilon}$, the structural Banach algebra of $W_{\epsilon}$, by a localization technique of holomorphic functional analysis. Now $\algnorm{\ndot}_{\phi, X|Y}$ is bounded from above by $\algnorm{\ndot}_{W_{\epsilon}}$ thanks to the continuity. Then we manage to include $W_{\epsilon}$ into $\M(\widehat{V}_{\sbullet}(L_{X|Y}, \phi|_Y(\epsilon)))$, so that $\algnorm{\ndot}_{W_{\epsilon}}$ is bounded by $\algnorm{\ndot}_{\phi|_Y(\epsilon)}$ from above, as they are both supremum norms on the corresponding domains. 

\subsection{Localization of spectrum by affinoid domain covering}\hfill\break
For any $\epsilon>0$, one can identify $\M(\widehat{V}_{\sbullet}(L_{X|Y}, \phi|_Y(\epsilon)))$ with $\M(\widehat{V}_{\sbullet}(L_{X|Y}, \phi(\epsilon)_{X|Y}))$ by Corollary \ref{spectral of quot is sup}. Via this identification, one has 
\[\M(\widehat{V}_{\sbullet}(L_{X|Y}, \phi_{X|Y}))\subseteq \M(\widehat{V}_{\sbullet}(L_{X|Y}, \phi|_Y(\epsilon))).\]
By Corollary \ref{general box all geom}, they can be both identified with closed subsets in $\M(\widehat{V}_{\sbullet}(L_{X|Y}, \phi(\epsilon)_{X|Y}^{\mathrm{aff}}))$, which itself can be identified with a closed subset in $\spec(V_{\sbullet}(L_{X|Y}))^{\mathrm{an}}$. We shall work in this fixed affinoid domain. We use notations $\Int^{\mathrm{top}}_{\M}$ and $\Int^{\mathrm{top}}_{V_{\sbullet}}$ to distinguish the topological interior in $\M(\widehat{V}_{\sbullet}(L_{X|Y}, \phi(\epsilon)_{X|Y}^{\mathrm{aff}}))$ and in $\spec(V_{\sbullet}(L_{X|Y}))^{\mathrm{an}}$.

\begin{lemma}\label{disc bundle interior}
    Assume that $\phi|_Y$ is asymptotic Fubini-Study. For any $\epsilon>0$,  $\M(\widehat{V}_{\sbullet}(L_{X|Y}, \phi_{X|Y}))$ is contained in $\Int^{\mathrm{top}}_{\M}(\M(\widehat{V}_{\sbullet}(L_{X|Y}, \phi|_Y(\epsilon))))$.
\end{lemma}
\begin{proof}
    By Proposition \ref{Envelop quot=Envelop sup}, $\mathcal{P}(\algnorm{\ndot}_{\phi|_Y})$ is equal to $\phi|_Y$, so it is continuous. Hence for any $\epsilon>0$, by Proposition \ref{disc inclusion}, one has
    \[\M(\widehat{V}_{\sbullet}(L_{X|Y}, \phi_{X|Y}))=\M(\widehat{V}_{\sbullet}(L_{X|Y}, \phi|_Y))\subseteq \Int^{\mathrm{top}}_{V_{\sbullet}}(\M(\widehat{V}_{\sbullet}(L_{X|Y}, \phi|_Y(\epsilon)))).\]
    By Proposition \ref{normed spectrum into non-normed spectrum}, the topology on $\M(\widehat{V}_{\sbullet}(L_{X|Y}, \phi(\epsilon)_{X|Y}^{\mathrm{aff}}))$ coincides the induced topology from $\spec(V_{\sbullet}(L|_Y))^{\mathrm{an}}$, hence the set $\Int^{\mathrm{top}}_{V_{\sbullet}}(\M(\widehat{V}_{\sbullet}(L_{X|Y}, \phi|_Y(\epsilon))))$ which is open in $\spec(V_{\sbullet}(L|_Y))^{\mathrm{an}}$ is also open in $\M(\widehat{V}_{\sbullet}(L_{X|Y}, \phi(\epsilon)_{X|Y}^{\mathrm{aff}}))$. Therefore this set is contained in $\Int^{\mathrm{top}}_{\M}(\M(\widehat{V}_{\sbullet}(L_{X|Y}, \phi|_Y(\epsilon))))$.
\end{proof}

\begin{proposition}\label{special domain neighbourhood}
    Assume that $\phi|_Y$ is asymptotic Fubini-Study. For any $\epsilon>0$, there exist a special domain $W_{\epsilon}$ such that 
    \[\M(\widehat{V}_{\sbullet}(L_{X|Y}, \phi_{X|Y}))\subseteq W_{\epsilon}\subseteq \M(\widehat{V}_{\sbullet}(L_{X|Y}, \phi|_Y(\epsilon))).\]
\end{proposition}
\begin{proof}
    By Corollary \ref{affinoid topological basis}, every point in $\M(\widehat{V}_{\sbullet}(L_{X|Y}, \phi(\epsilon)_{X|Y}^{\mathrm{aff}}))$ has a neighbourhood system consisting of affinoid domains. Hence for any $z\in \M(\widehat{V}_{\sbullet}(L_{X|Y}, \phi_{X|Y}))$, there exists an affinoid domain neighbourhood $V(z)$. By Lemma \ref{disc bundle interior}, $z$ has an open neighbourhood $\Int^{\mathrm{top}}_{\M}(\M(\widehat{V}_{\sbullet}(L_{X|Y}, \phi|_Y(\epsilon))))$, so we can assume that each $V(z)$ is contained in this open set. 
    
    One forms a covering by open sets
    \[\M(\widehat{V}_{\sbullet}(L_{X|Y}, \phi_{X|Y}))\subseteq \bigcup_{z\in \M(\widehat{V}_{\sbullet}(L_{X|Y}, \phi_{X|Y}))}\Int^{\mathrm{top}}_{\M}V(z).\]
    Since the left hand side is a compact set by Proposition \ref{banach spectrum compact Hausdorff}, there exist finitely many points $\{z_1, \dots, z_m\}$ such that $\{\Int^{\mathrm{top}}_{\M}V(z_i)\}_{i\in \card{m}}$ form a covering of $\M(\widehat{V}_{\sbullet}(L_{X|Y}, \phi_{X|Y}))$. Let $W_{\epsilon}$ be the union of affinoid domains $\{V(z_i)\}_{i\in \card{m}}$, then it is a special domain, and satisfies the desired inclusion conditions.
\end{proof}
\begin{remark}
    One denotes by $\mathcal{W}_{\epsilon}$ the Banach $k$-algebra of $\Gamma(W_{\epsilon}, \mathscr{O}_{W_{\epsilon}})$ equipped with supremum norm $\algnorm{\ndot}_{W_{\epsilon}}$ (see Definition \ref{special domain}). 
\end{remark}
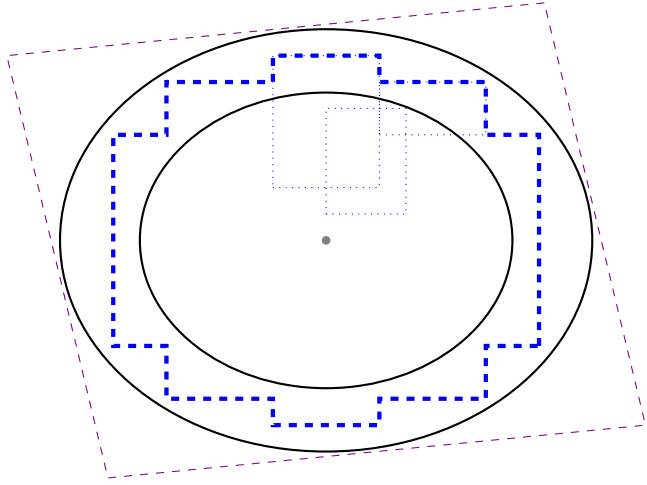
\begin{tikzpicture}[scale=0.70]

\filldraw [gray] (0,0) circle (2pt);
\draw [thick] (0, 0) ellipse (3.5 and 2.8);
\draw [thick] (0, 0) ellipse (5 and 4);

\draw [blue, dashed, ultra thick] (4,-2) -- (4,2) -- (3,2) -- (3,3) -- (1,3)-- (1, 3.5) -- (-1, 3.5) -- (-1,3) -- (-3,3) -- (-3,2) -- (-4,2) -- (-4,-2) -- (-3,-2) -- (-3,-3) -- (-1,-3) -- (-1, -3.5) -- (1, -3.5) -- (1, -3) -- (3, -3) -- (3,-2) -- cycle;

\draw [violet, dashed] (4.1, 4.5) -- (-6,3.5) -- (-4.1,-4.5) -- (6,-3.5) -- cycle;

\draw [blue, dotted] (1,3) rectangle (3,2);

\draw [blue, dotted] (-1,3.5) rectangle (1,1);
\draw [blue, dotted] (0,2.5) rectangle (1.5,0.5);

\node [below] at (0, -5) {\textbf{Geometric approximation}};
\node [below] at (0,-6)
{$\M(\widehat{V}_{\sbullet}(L_{X|Y}, \phi_{X|Y}))\subseteq {\color{blue}W_{\epsilon}}\subseteq \M(\widehat{V}_{\sbullet}(L_{X|Y}, \phi|_{Y}(\epsilon)))\subseteq {\color{violet}
\M(\widehat{V}_{\sbullet}(L_{X|Y}, \phi_{X|Y}(\epsilon)^{\mathrm{aff}})}$};

\end{tikzpicture}

\subsection{Localization of Banach algebra homomorphism}\hfill\break
For any $\epsilon>0$, by Theorem \ref{special domain neighbourhood} for $\epsilon/2$, there exist a special domain $W_{\epsilon/2}$ such that
\[\M(\widehat{V}_{\sbullet}(L_{X|Y}, \phi_{X|Y}))\subseteq W_{\epsilon/2}\subseteq \M(\widehat{V}_{\sbullet}(L_{X|Y}, \phi|_Y(\epsilon/2))).\]

\begin{proposition}\label{QE spectral localization}
    For any $\epsilon>0$, there exist a homomorphism of Banach $k$-algebras
    \[\theta(\epsilon/2)_{W_{\epsilon/2}}: \mathcal{W}_{\epsilon/2}\rightarrow \widehat{V}_{\sbullet}(L_{X|Y}, \phi_{X|Y})\]
    which extends the identity map on the dense sub-$k$-algebra $V_{\sbullet}(L_{X|Y})$.
\end{proposition}
\begin{proof}
    By Proposition \ref{general box all alg}, there are homomorphisms of Banach algebras induced by the identity map on the dense $V_{\sbullet}(L_{X|Y})$:
    \[\widehat{V}_{\sbullet}(L_{X|Y}, \phi(\epsilon/2)_{X|Y}^{\mathrm{aff}})\xrightarrow[\sigma(\epsilon/2)|_Y]{} \widehat{V}_{\sbullet}(L_{X|Y}, \phi(\epsilon/2)_{X|Y})\xrightarrow[\iota(\epsilon/2)]{} \widehat{V}_{\sbullet}(L_{X|Y}, \phi_{X|Y}).\]
    Let $\tau(\epsilon/2)$ denote the composed homomorphism of Banach $k$-algebras. It is a homomorphism from an affinoid algebra to a Banach algebra. It has dense image, so by Proposition \ref{injective} the induced continuous map 
    \[\tau(\epsilon/2)^*: \M(\widehat{V}_{\sbullet}(L_{X|Y}, \phi_{X|Y}))\rightarrow \M(\widehat{V}_{\sbullet}(L_{X|Y}, \phi(\epsilon/2)_{X|Y}^{\mathrm{aff}}))\] 
    is injective and is closed. As both spaces are compact and Hausdorff, this map is a homeomorphism from its domain to its image. So the homomorphism spectrum $\Sigma_{\tau(\epsilon/2)}$ is homeomorphic to $\M(\widehat{V}_{\sbullet}(L_{X|Y}, \phi_{X|Y}))$, and is contained in $W_{\epsilon/2}$.
    
    One performs spectral calculus for the homomorphism $\tau(\epsilon/2)$ and the special domain $W_{\epsilon/2}$: by Theorem \ref{holomorphic functional calculus}, there exist a homomorphism of Banach $k$-algebras
    \[\theta(\epsilon/2)_{W_{\epsilon/2}}: \mathcal{W}_{\epsilon/2}\rightarrow \widehat{V}_{\sbullet}(L_{X|Y}, \phi_{X|Y})\]
    which extends the identity map on the dense sub-$k$-algebra $V_{\sbullet}(L_{X|Y})$.
\end{proof}

\begin{theorem}\label{quantitative extension geometric}
    Let $\phi$ be an asymptotic Fubini-Study metric on $L$. Then for any $\epsilon>0$, there exists $n_Y\in\N$ such that for any $n\geq n_Y$ and any $t_n\in V_n(L|_Y)$, there exits $s_n\in V_n(L)$ such that $s_n|_Y=t_n$ and 
    \[\norm{s_n}_{n\phi}\leq \mathrm{e}^{n\epsilon}\cdot\norm{t_n}_{n\phi|_Y}.\]
\end{theorem}
\begin{proof}
    Start from the homomorphism $\theta(\epsilon/2)_{W_{\epsilon/2}}$ constructed in Proposition \ref{QE spectral localization}. The boundedness (continuity) of this homomorphism of Banach $k$-algebras implies that there exists $C_{\epsilon}>0$ such that
    \[\forall \underline{t}\in V_{\sbullet}(L_{X|Y}), \quad  \algnorm{\underline{t}}_{\phi, X|Y}\leq C_{\epsilon}\cdot \algnorm{\underline{t}}_{W_{\epsilon/2}}=C_{\epsilon}\cdot \sup_{z\in W_{\epsilon/2}\cdot}\abs{\underline{t}}_z.\]
    
    By Proposition \ref{Envelop quot=Envelop sup}, one has a canonical homeomorphism
    \[\M(\widehat{V}_{\sbullet}(L_{X|Y}, \phi(\epsilon/2)_{X|Y}))\simeq \M(\widehat{V}_{\sbullet}(L_{X|Y}, \phi|_Y(\epsilon/2)))\]
    from which one deduces 
    \[W_{\epsilon/2}\subseteq \M(\widehat{V}_{\sbullet}(L_{X|Y}, \phi|_Y(\epsilon/2))).\]
    Remember that since $\algnorm{\ndot}_{\phi|_Y(\epsilon/2)}$ is power-mutliplicative, by Theorem \ref{Gelfand uniform}, it is the supremum norm on $\M(\widehat{V}_{\sbullet}(L_{X|Y}, \phi|_Y(\epsilon/2)))$. Hence by comparing supremum norms on these two closed sets, we get
    \[\forall \underline{t}\in V_{\sbullet}(L_{X|Y}), \quad \algnorm{\underline{t}}_{W_{\epsilon/2}}\leq \algnorm{\underline{t}}_{\phi|_Y(\epsilon/2)},\]
    therefore for any $t_n \in V_{n}(L_{X|Y})$, one has
    \begin{equation*}
        \begin{split}
            \norm{t_n}_{n\phi, X|Y}&= \algnorm{t_n}_{\phi, X|Y}\\
            &\leq C_{\epsilon}\cdot \algnorm{t_n}_{W_{\epsilon/2}}\\
            &\leq C_{\epsilon}\cdot \algnorm{t_n}_{\phi|_Y(\epsilon/2)}=C_{\epsilon}\cdot\mathrm{e}^{n\epsilon/2}\cdot\norm{t_n}_{n\phi|_Y}.
        \end{split}
    \end{equation*}
    Let $n_Y$ be the integer $\max\{\lceil\log(C_{\epsilon})/(\epsilon/2)\rceil, M\}$, then for any $n\geq n_Y$ and any $t_n \in V_{n}(L|_Y)$, there exists $s_n\in V_n(L)$ such that
    \[\norm{s_n}_{n\phi}\leq \mathrm{e}^{n\epsilon}\cdot \norm{t_n}_{n\phi|_Y}.\]
\end{proof}

\section{Algebraic approximation}

In this section, we approximate the Banach algebra norm $\algnorm{\ndot}_{\phi_{X|Y}}$ by affinoid norms. This algebraic approximation exploits subtle consequences of existence of a ultra-metric orthogonal basis in $(V_n(L), \norm{\ndot}_{n\phi})$ for some $n\in \N$.

\emph{In this section $(k, \abs{\ndot})$ is assumed to be discretely valued}. With this assumption, recall that $d+1$ classes of real numbers $\parenth{\alpha(p_0), \dots, \alpha(p_d)}$ in $\R/H(k, \abs{\ndot})$ are said to be $\Q$-independent if there exists no $(a_0,\dots, a_d)\in \Q^{d+1}$ and no $p\in H(k, \abs{\ndot})$ such that $\sum_{i\in \Card{d}}a_i\cdot p_i=p$. Recall that thanks to the discreteness of $\abs{\ndot}$, by Proposition \ref{discrete orthogonal vector space}, any finite dimensional ultrametric normed $k$-vector space has an orthogonal basis.

\subsection{Algebra norm induced by Fubini-Study metric}

\subsubsection{Case for $(\mathbb{P}^d, \mathscr{O}(1))$}\hfill\break

Let $\phi$ be a metric on $\mathscr{O}(1)$, one studies the algebra norm $\algnorm{\ndot}_{\phi}$ on $V_{\sbullet}(\mathscr{O}(1))$. One would like to show that with various assumptions, it is a Gauss algebra norm, namely the standard affinoid algebra norm on the polynomial algebra. Then the normed section algebra $(V_{\sbullet}(\mathscr{O}(1)), \algnorm{\ndot}_{\phi})$ will be a Tate affinoid algebra. (see Definition \ref{standard Tate})

By Proposition \ref{discrete orthogonal vector space}, there exist an orthogonal basis $\{T_i\}_{i\in \Card{d}}$ for the normed vector space $(V_1(\mathscr{O}(1)), \norm{\ndot}_{\phi})$. For any $i\in \Card{d}$, one denotes by $r_i$ the value $\norm{T_i}_{\phi}$, and by $\pmb{r}\in (\R_+)^{d+1}$ the multi-radius $(r_0, \dots, r_d)$. One fixes such an orthogonal basis, and identify the graded $k$-algebra $V_{\sbullet}(\mathscr{O}(1))$ with $k[T_0,\dots, T_d]$. For any multi-index $J=(j_0, \dots, j_d)\in \N^{d+1}$, one denotes by $\pmb{T}^{J}$ the monomial element $\prod_{i\in \Card{d}}(T_i)^{j_i}\in V_{\abs{J}}(\mathscr{O}(1))$. 

Note that in general, the sub-spaces $V_n(\mathscr{O}(1))$ are orthogonal with respect to $\algnorm{\ndot}_{\phi}$ for different $n\in \N$, while a Gauss algebra norm exhibits a much finer orthogonality: the sub-spaces generated by each mononial $\pmb{T}^J$ should be orthogonal for different $J\in \N^{d+1}$.

First, on monomial elements, the algebra norm $\algnorm{\ndot}_{\phi}$ resembles a Gauss norm.
\begin{proposition}\label{monomial simultaneous maximum}
    For any $J\in \N^{d+1}$, one has
    \[\norm{\pmb{T}^{J}}_{|J|\phi}=\prod_{i\in \Card{d}}\norm{T_i}_{\phi}^{j_i}.\]
\end{proposition}
\begin{proof}
    Take a complete non-Archimedean valued field extension $(K, \abs{\ndot}_K)$ of $(k, \abs{\ndot}_k)$ such that
    \[\forall i\in \Card{d}, \quad \{r_i\}_{i\in \Card{d}}\subseteq\abs{k^{\times}}_K,\]
    hence for any $i\in \Card{d}$, there exist elements $\kappa_i\in K$ such that $\abs{\kappa_i}_K=r_i$. One denotes by $x(\pmb{r})\in (\mathbb{P}^d_k)^{\mathrm{an}}$ the point given by coordinates $[\kappa_0:\dots:\kappa_d]$.
    
    \begin{claim}
        For any $i\in \Card{d}$, one has
        \[\norm{T_i}_{\phi}=r_i=\abs{T_i}_{\phi}(x(\pmb{r})).\]
        In other words, the maximum of the function $\abs{T_i}_{\phi}(x)$ on $(\mathbb{P}^d_k)^{\mathrm{an}}$ is $r_i$, and the maximum values of these $d+1$ functions can be attained \textbf{at the same point} $x(\pmb{r})$.
    \end{claim}
    \begin{proof}
    By the orthogonality of the basis $\{T_i\}_{i\in \Card{d}}$, we can compute
        \begin{equation*}
        \begin{split}
            \abs{T_i}_{\phi}(x(\pmb{r}))&=\inf_{\begin{subarray}{c}(\sum_{m\in \Card{d}}f_m\cdot T_m)(x(\pmb{r}))=(T_i)(x(\pmb{r}))\\
            (f_0, \dots, f_d)\in k^{d+1}
            \end{subarray}}\Big\lVert\sum_{m\in \Card{d}}f_m\cdot T_m\Big\rVert_{\phi}
            \\&= \inf\limits_{\sum_{m\in \Card{d}}f_m\kappa_m=\kappa_i}\max_{m\in \Card{d}}\Big\{\norm{f_m\cdot T_m}_{\phi}\Big\}
            \\&= \inf\limits_{\sum_{m\in \Card{d}}f_m\kappa_m=\kappa_i}\max_{m\in \Card{d}}\Big\{\abs{f_m}\abs{\kappa_m} \Big\}
            \\&= \inf\limits_{\sum_{m\in \Card{d}}f_m(\kappa_m/\kappa_i)=1}\max_{m\in \Card{d}}\Big\{\abs{f_m}\abs{\kappa_m/\kappa_i}\cdot \abs{\kappa_i} \Big\}
            \\&=\abs{\kappa_i}=r_i.
        \end{split}
        \end{equation*}
    The last equality is obtained by Lemma \ref{ultrametric ortho min max}. 
    \end{proof}
    
    By this Claim, for any multi-index $J$, the function $\abs{\pmb{T}^{J}}_{\phi}(x)$ can attain its maximum value $\prod_{i\in \Card{d}}r_i^{j_i}$ at the point $x(\pmb{r})\in(\mathbb{P}_k^d)^{an}$ as the product of maximum of factors of the monomial. By definition, 
    \[\norm{\pmb{T}^{J}}_{|J|\phi}=\sup_{x\in (\mathbb{P}^d)^{an}}\abs{\pmb{T}^{J}}_{\phi}(x)=\prod_{i\in \Card{d}}r_i^{j_i}=\prod_{i\in \Card{d}}\norm{T_i}_{\phi}^{j_i}.\]
\end{proof}

Second, one calculates the algebra norm $\algnorm{\ndot}_{\phi}$ on any (homogeneous) combination of monomials. For a general metric, one needs a $\Q$-independence assumption to gain finer orthogonality.
\begin{proposition}\label{monomial orthogonal}
    Assume that $\{\alpha(\norm{T_i}_{\phi})\}_{i\in \Card{d}}$ are $\Q$-independent in $\R/H(k, \abs{\ndot})$. Let $S\subseteq \N^{d+1}$ be a finite set of multi-indices, then for any $J\in S$ any $f_{J}\in k$, one has
    \[\Big\vvvert\sum_{\begin{subarray}{c}J\in S \end{subarray}}f_{J}\cdot \pmb{T}^{J}\Big\vvvert = \sup_{\begin{subarray}{c}J\in S \end{subarray}}\  \norm{f_{J}\cdot \pmb{T}^{J}}_{|J|\phi}.\]
    In other words, the algebra norm $\algnorm{\ndot}_{\phi}$ on $V_{\sbullet}(\mathscr{O}(1))$ is a Gauss norm on $k[T_0,\dots, T_d]$ of multi-radius $\pmb{r}$. The Banach $k$-algebra $\widehat{V}_{\sbullet}(\mathscr{O}(1), \phi)$ is an affinoid algebra.
\end{proposition}
\begin{proof}
    By the $\Q$-independence assumption and Proposition \ref{monomial simultaneous maximum}, for any two distinct multi-index $J$ and $J'$, and any two non-zero coefficients $f_{J}$ and $f_{J'}$ in $k$, we have
    \[\norm{f_{J}\cdot\pmb{T}^{J}}_{|J|\phi}=\abs{f_{J}}\cdot\prod_{i\in \Card{d}}r_i^{j_i}\neq \abs{f_{J'}}\cdot \prod_{i\in \Card{d}}r_i^{j'_i}=\norm{f_{J'}\cdot\pmb{T}^{J'}}_{|J'|\phi}.\]
    By Lemma \ref{distinct orthogonal}, the elements $\{\pmb{T}^{J}\}_{J\in S}$ form an orthogonal basis for the normed vector space $(\bigoplus_{J\in S}k\cdot \pmb{T}^{J}, \algnorm{\ndot}_{\phi})$. So the equality in the conclusion holds.
\end{proof}

\begin{corollary}\label{metric Q-indep on P^d is FS}
    With the same assumptions as above, the envelop metric $\mathcal{P}(\phi)$ is a Fubini-Study metric induced by $\norm{\ndot}_{\phi}$, and is continuous.
\end{corollary}

For a Fubini-Study metric, one does not need the $\Q$-independence. For any $\pmb{\delta}=(\delta_0,\dots, \delta_d)\in \R^{d+1}$, one constructs a perturbed metric $\phi(\pmb{\delta})$ as follows. Let $\norm{\ndot}_{\phi(\pmb{\delta})}$ be the norm on $V_1(\mathscr{O}(1))$ such that $\{T_i\}_{i\in \card{d}}$ is an orthogonal basis with new norms
\[\forall i\in \Card{d},\quad \norm{T_i}_{\phi}(\pmb{\delta})=\mathrm{e}^{\delta_i}\norm{T_i}_{\phi}.\]
Let $\phi(\pmb{\delta})$ be the metric $\mathrm{FS}(\norm{\ndot}_{\phi(\pmb{\delta})})$ on $\mathscr{O}(1)$. Let $\abs{\pmb{\delta}}$ denote the number $\max_{i\in \Card{d}}\abs{\delta_i}\in \R_+$.

\begin{lemma}\label{approximation Tate}
    Assume that $\phi$ is a Fubini-Study metric. For any $\epsilon>0$, there exists $\pmb{\delta}\in \R^{d+1}$ with $|\pmb{\delta}|\leq \epsilon$ such that
    \[\forall n\in \N,\quad \dist (\norm{\ndot}_{n\phi}, \norm{\ndot}_{n\phi(\pmb{\delta})})\leq n\epsilon.\]
\end{lemma}
\begin{proof}
    Choose an arbitrary $\pmb{\delta}$ with $|\pmb{\delta}|\leq \epsilon$. Then
    \[\dist(\norm{\ndot}_{\phi}, \norm{\ndot}_{\phi(\pmb{\delta})})\leq \epsilon. \]
    By Proposition \ref{FS envelop idempotent}, we have
    \[\dist(\mathrm{FS}(\norm{\ndot}_{\phi}), \mathrm{FS}(\norm{\ndot}_{\phi(\pmb{\delta})}))=\dist(\mathrm{FS}(\norm{\ndot}_{\phi}), \phi(\pmb{\delta}))\leq \epsilon,\]
    by the assumption and Proposition \ref{FS idempotent}, \[\mathrm{FS}(\norm{\ndot}_{\phi})=\phi,\]
    so the conlusion holds.
\end{proof}

\begin{proposition}\label{FS is Tate}
    Assume that $\phi$ is a Fubini-Study metric. The conclusion of Proposition \ref{monomial orthogonal} holds without the assumption of $\Q$-independence of $\{\alpha(\norm{T_i}_1)\}_{i\in \Card{d}}$.
\end{proposition}
\begin{proof}
    Since $\abs{\ndot}_k$ is discrete, for any $\epsilon>0$, there exists $\pmb{\delta}$ with $|\pmb{\delta}|\leq \epsilon$ such that the elements $\{\alpha(\norm{T_i}_{\phi(\pmb{\delta})})\}_{i\in \Card{d}}$ are $\Q$-independent in $\R/H(k, \abs{\ndot})$. By Proposition \ref{approximation Tate}, for any $n\in \N$ and any $s_n=\sum_{\abs{J}=n}f_{J}\cdot \pmb{T}^{J}\in V_n(\mathscr{O}(1))$,
    \[\mathrm{e}^{-n\epsilon}\Big\lVert\sum_{\abs{J}=n}f_{J}\cdot \pmb{T}^{J}\Big\rVert_{n\phi(\pmb{\delta})}\leq \Big\lVert\sum_{\abs{J}=n}f_{J}\cdot \pmb{T}^{J}\Big\rVert_{n\phi} \leq \mathrm{e}^{n\epsilon}\Big\lVert\sum_{\abs{J}=n}f_{J}\cdot \pmb{T}^{J}\Big\rVert_{n\phi(\pmb{\delta})}.\]
    By Proposition \ref{monomial orthogonal}, one has
    \[\max_{\abs{J}=n}\Big\{\mathrm{e}^{-n\epsilon}\abs{f_{J}}\prod_{i\in \Card{d}}(\mathrm{e}^{\delta_i}r_i)^{j_i}\Big\}
    \leq \Big\lVert\sum_{\abs{J}=n}f_{J}\cdot \pmb{T}^{J}\Big\rVert_{n\phi}
    \leq \max_{\abs{J}=n}\Big\{\mathrm{e}^{n\epsilon}\abs{f_{J}}\prod_{i\in \Card{d}}(\mathrm{e}^{\delta_i}r_i)^{j_i}\Big\}.\]
    Fix $n$ and let $\epsilon\to 0$, one gets
    \[\Big\lVert\sum_{\abs{J}=n}f_{J}\cdot \pmb{T}^{J}\Big\rVert_{n\phi}=\max_{\abs{J}=n}\Big\{\abs{f_{J}}\prod_{i\in \Card{d}}r_i^{j_i}\Big\}.\]
    So one has
    \[\Big\vvvert\sum_{\abs{J}<\infty}f_{J}\cdot \pmb{T}^{J}\Big\vvvert_{n\phi}=\sup_{n\in \N}\max_{\abs{J}=n}\Big\{\abs{f_{J}}\prod_{i\in \Card{d}}r_i^{j_i}\Big\}=\max_{\abs{J}<\infty}\Big\{\abs{f_{J}}\prod_{i\in \Card{d}}r_i^{j_i}\Big\}.\]
    Hence $\algnorm{\ndot}_{\phi}$ is a Gauss norm of multi-radius $\pmb{r}$ on $V_{\sbullet}(\mathscr{O}(1))$.
\end{proof}

\subsubsection{Case for general $(X, L)$}

\begin{proposition}\label{very ample FS is affinoid}
    Assume that $\phi$ is a Fubini-Study metric. If $L$ is very ample, then $\widehat{V}_{\sbullet}(L, \phi)$, $\widehat{V}(L_{X|Y}, \phi_{X|Y})$ and $\widehat{V}_{\sbullet}(L|_Y, \phi|_Y)$ are affinoid algebras. 
\end{proposition}
\begin{proof}
    By the assumption, the elements of $V_1(L)$ induces an embedding
    \[\iota_1: X\rightarrow \mathbb{P}^{d_1}_k\]
    such that $\iota_1^*\mathscr{O}(1)=L$ with $\mathrm{dim}_k V_1=d_1+1$. Moreover there exists a norm $\norm{\ndot}_1$ on $V_1(L)$ such that $\phi=\mathrm{FS}(\norm{\ndot}_1)$. View $\norm{\ndot}_1$ as a norm on $V_1(\mathscr{O}(1))$, we get a metric $\psi=\mathrm{FS}(\norm{\ndot}_1)$ on $\mathscr{O}(1)$. By construction $\psi|_X=\phi$.
    
    By Proposition \ref{FS is Tate}, the Banach algebra $\widehat{V}_{\sbullet}(\mathscr{O}(1), \algnorm{\ndot}_{\psi})$ is an affinoid algebra. Hence the quotient Banach algebra $\widehat{V}_{\sbullet}(L, \algnorm{\ndot}_{\psi,\mathbb{P}^{d_1}|X})$ is an affinoid algebra.
    
    By Proposition \ref{spectral of quot is sup}, the algebra norm $\algnorm{\ndot}_{\phi}$ is the spectral norm of $\algnorm{\ndot}_{\psi,\mathbb{P}^{d_1}|X}$. The later is an affinoid norm, hence is equivalent to its spectral norm by Proposition \ref{affinoid norm equivalent to spectral norm}. So $\algnorm{\ndot}_{\phi}$ is also an affinoid algebra norm. Thus $\widehat{V}_{\sbullet}(L, \phi)$ is an affinoid algebra.
    
    Similarly, by Proposition \ref{spectral of quot is sup}, on $V_{\sbullet}(L_{X|Y})$, the algebra norm $\algnorm{\ndot}_{\phi|_Y}$ is the spectral norm of $\algnorm{\ndot}_{\phi, X|Y}$, hence is itself an affinoid algebra norm.
\end{proof}

\begin{corollary}\label{ample FS is affinoid}
    Assume that $\phi$ is a Fubini-Study metric. If $L$ is just ample, then $\widehat{V}_{\sbullet}(L, \phi)$, $\widehat{V}(L_{X|Y}, \phi_{X|Y})$ and $\widehat{V}_{\sbullet}(L|_Y, \phi|_Y)$ are affinoid algebras.
\end{corollary}
\begin{proof}
    By assumption, $L^{\otimes M}$ is very ample. So $\widehat{V}_{\sbullet}^{(M)}(L, \phi)$, $\widehat{V}_{\sbullet}^{(M)}(L_{X|Y}, \phi_{X|Y})$ and $\widehat{V}_{\sbullet}^{(M)}(L|_Y, \phi|_Y)$ are affinoid algebras. Since $V_{\sbullet}(L)$ is integral and is finite over $V_{\sbullet}^{(M)}(L)$, by Proposition \ref{Affinoid sub-algebra dense finite is finite}, the Banach algebras $\widehat{V}_{\sbullet}(L, \phi)$, $\widehat{V}(L_{X|Y}, \phi_{X|Y})$ and $\widehat{V}_{\sbullet}(L|_Y, \phi|_Y)$ are Banach finite over $\widehat{V}_{\sbullet}^{(M)}(L, \phi)$, $\widehat{V}_{\sbullet}^{(M)}(L_{X|Y}, \phi_{X|Y})$ and $\widehat{V}_{\sbullet}^{(M)}(L|_Y, \phi|_Y)$ respectively. Hence they are affinoid algebras.  
\end{proof}

\begin{proposition}\label{QE algebraic FS}
    Assume that $\phi$ is a Fubini-Study metric. Then there exist $C(\phi, Y)>0$ such that for any $\underline{t}\in V_{\sbullet}(L_{X|Y})$, there exists $\underline{s}\in V_{\sbullet}(L)$ with
    \[\algnorm{\underline{s}}_{\phi}\leq C(\phi,Y, X)\cdot \algnorm{\underline{t}}_{\phi|_Y}.\]
    In particular, for any $n\in \N$ and $t_n\in V_{\sbullet}(L_{X|Y})$, there exists $s_n\in V_{\sbullet}(L)$ with
    \[\norm{s_n}_{n\phi}\leq C(\phi, Y, X)\cdot \norm{t_n}_{n\phi|_Y}.\]
\end{proposition}
\begin{proof}
    By Corollary \ref{ample FS is affinoid}, the Banach algebra norm $\algnorm{\ndot}_{\phi_{X|Y}}$ is an affinoid algebra norm. By Proposition \ref{affinoid norm equivalent to spectral norm}, there exists $C(\phi, Y)>0$ such that
    \[\algnorm{\ndot}_{\phi_{X|Y}}\leq C(\phi, Y, X)\cdot \algnorm{\ndot}_{\phi_{X|Y}; \mathrm{sp}}.\]
    Since $\algnorm{\ndot}_{\phi_{X|Y}; \mathrm{sp}}=\algnorm{\ndot}_{\phi|_Y}$ by Corollary \ref{spectral of quot is sup}, one gets the bounds.
\end{proof}
\begin{remark}\label{affinoid estimate vs OT L2}
    With the metric finiteness properties of affinoid algebra norm, here the upper bound for metric extension of a Fubini-Study metric is much better than what was expected, compared to (\ref{Equ: Bost-Randriam bound}) or even to (\ref{Equ: Tian bound}), for its (in)depence on $n\in \N$. This independence suggest that it would be reasonable to compare this affinoid algebra technique in this non-Archimedean setting with the use of Ohsawa-Takegoshi $L^2$ extension technique in the complex analytic setting.
\end{remark}

\subsection{Algebra norm induced by asymptotic Fubini-Study metric}\hfill\break
With the extra assumption of discreteness for the base valued field, one can give another proof of Theorem \ref{quantitative extension geometric}.

\begin{theorem}\label{quantitative extension algebraic}
    Suppose that $(k, \abs{\ndot})$ is discretely valued. Let $\phi$ be an asymptotic Fubini-Study metric on $L$. Then for any $\epsilon>0$, there exists $n_Y\in\N$ such that for any $n\geq n_Y$ and any $t_n\in V_n(L|_Y)$, there exits $s_n\in V_n(L)$ such that $s_n|_Y=t_n$ and 
    \[\norm{s_n}_{n\phi}\leq \mathrm{e}^{n\epsilon}\cdot\norm{t_n}_{n\phi|_Y}.\]
\end{theorem}
\begin{proof}
    Recall that one can find $M\in \N$ such that $L^{\otimes M}$ is very ample and for any $n\geq M$, the restriction map from $V_n(L)$ to $V_n(L|_Y)$ is surjective.
    
    By the asymptotic Fubini-Study assumption, there exist norms $\{\norm{\ndot}_n\}_{n\in \N}$ on $V_n(L)$ such that $\phi$ is given by $\mathcal{P}(\{\norm{\ndot}_n\}_{n\in \N})$. Let $\psi_n$ denote the metric $\frac{1}{n}\mathrm{FS}(\norm{\ndot})_n$ on $L$. By the continuity assumption, the convergence to envelop metric is uniform (see \S3.1 \ref{Item: FS metric and FS envelop metric}). So for any $\epsilon>0$, there exists $M'\geq M$ such that
    \[\dist (\psi_{M'}, \phi)\leq \frac{1}{3}\epsilon, \]
    hence for any $n\in \N$, one has
    \[\dist (\norm{\ndot}_{n\psi_{M'}},  \norm{\ndot}_{n\phi})\leq \frac{1}{3}n\epsilon, \quad \dist (\norm{\ndot}_{n\psi_{M'}, X|Y},  \norm{\ndot}_{n\phi, X|Y})\leq \frac{1}{3}n\epsilon. \]
    By Corollary \ref{ample FS is affinoid}, the Banach algebra $\widehat{V}_{\sbullet}(L_{X|Y}, (\psi_{M'})_{X|Y})$ is an affinoid algebra. By Proposition \ref{QE algebraic FS}, there exists $C_{M'}(X, Y, \phi)>0$ such that
    \[\forall \underline{t}\in V_{\sbullet}(L_{X|Y}), \quad \algnorm{\underline{t}}_{\psi_{M'}, X|Y}\leq C_{M'}\cdot \algnorm{\underline{t}}_{\psi_{M'}, X|Y;\mathrm{sp}}=\algnorm{\underline{t}}_{\psi_{M'}, Y}.\]
    Combining this comparison with above estimates, one has that for every $n\geq n_Y:=\lceil \ln(C_{M'})/(\epsilon/3)\rceil$
    \begin{equation*}
        \begin{split}
            \norm{t_n}_{n\phi, X|Y}&\leq \mathrm{e}^{\frac{1}{3}n\epsilon}\cdot \norm{t_n}_{n\psi_{M'}, X|Y}
            \\&\leq \mathrm{e}^{\frac{1}{3}n\epsilon}\cdot C_{M'}\cdot  \norm{t_n}_{n\psi_{M'}|_Y}
            \\&\leq \mathrm{e}^{\frac{2}{3}n\epsilon}\cdot C_{M'}\cdot  \norm{t_n}_{n\phi|_Y}
            \\&\leq \mathrm{e}^{n\epsilon}\cdot  \norm{t_n}_{n\phi|_Y}.
        \end{split}
    \end{equation*}
    Hence there exists $s_n\in V_n(L)$ with $\norm{s_n}_{n\phi}\leq \mathrm{e}^{n\epsilon}\cdot  \norm{t_n}_{n\phi|_Y}$.
\end{proof}

\medskip


\end{document}